\documentclass[12pt]{amsart}
\usepackage{amsmath,amscd,amssymb,amsfonts,latexsym}
\usepackage[all]{xy}
\usepackage{ulem}

\usepackage{geometry}
\allowdisplaybreaks

\newcommand{\bC}{{\mathbb C}}

\newcommand{\bP}{{\mathbb P}}
\newcommand{\bR}{{\mathbb R}}
\newcommand{\bQ}{{\mathbb Q}}
\newcommand{\bZ}{{\mathbb Z}}

\newcommand{\cA}{{\mathcal A}}

\newcommand{\cE}{{\mathcal E}}
\newcommand{\cF}{{\mathcal F}}

\newcommand{\cH}{{\mathcal H}}

\newcommand{\cO}{{\mathcal O}}

\newcommand{\cS}{{\mathcal S}}

\newcommand{\BM}{{\rm BM}}

\newcommand{\Hom}{{\rm{Hom}}}
\newcommand{\Sp}{{\rm{Spec}}}

\newcommand{\wti}{\widetilde}

\newcommand{\Gr}{\text{\rm Gr}}

\newcommand{\lra}{\longrightarrow}

\newcommand{\con}{{\it Cone}}
\newcommand{\Int}{{\rm Int}}
\newcommand{\Relint}{{\rm Relint}}

\newcommand{\q}{\quad}

\newcommand{\ch}{{\rm ch}}
\newcommand{\Td}{{\rm Td}}

\newcommand{\mult}{{\rm mult}}

\newcommand{\os}{$O_{\sigma}$}
\newcommand{\vs}{$V_{\sigma}$}

\newcommand{\tyh}{$\widehat{T}_{y*}$}

\newcommand{\om}{\Omega}

\newcommand{\La}{\Lambda}

\def\al{\alpha}

\def\sig{\sigma}
\def\Sig{\Sigma}

\theoremstyle{plain}
\newtheorem{thm}{Theorem}[section]
\newtheorem{cor}[thm]{Corollary}

\newtheorem{lem}[thm]{Lemma}
\newtheorem{prop}[thm]{Proposition}

\newtheorem*{ack}{Acknowledgements}

\theoremstyle{definition}
\newtheorem{df}[thm]{Definition}
\newtheorem{rem}[thm]{Remark}
\newtheorem{example}[thm]{Example}

\def\be{\begin{equation}}
\def\ee{\end{equation}}

\def\bt{\begin{thm}}
\def\et{\end{thm}}

\def\bc{\begin{cor}}
\def\ec{\end{cor}}

\def\br{\begin{rem}}
\def\er{\end{rem}}

\def\bp{\begin{prop}}
\def\ep{\end{prop}}

\def\bl{\begin{lem}}
\def\el{\end{lem}}

\def\bn{\begin{enumerate}}
\def\en{\end{enumerate}}

\def\bex{\begin{example}}
\def\eex{\end{example}}

\def\bd{\begin{df}}
\def\ed{\end{df}}

\numberwithin{equation}{section}

\title{Characteristic classes of singular toric varieties}
\author{Lauren\c{t}iu Maxim}
\address{L. Maxim: Department of Mathematics, University of Wisconsin, 480 Lincoln Drive, Madison, WI 53706, USA}
\email {maxim@math.wisc.edu}
\author[J. Sch\"urmann ]{J\"org Sch\"urmann}
\address{J.  Sch\"urmann : Mathematische Institut,
          Universit\"at M\"unster,
          Einsteinstr. 62, 48149 M\"unster,
          Germany.}
\email {jschuerm@math.uni-muenster.de}

\date{\today}



\begin{document}

\begin{abstract}
In this paper we compute the motivic Chern classes and homology Hirzebruch characteristic classes of (possibly singular) toric varieties, which in the context of complete toric varieties fit 
nicely with a generalized Hirzebruch-Riemann-Roch theorem. As important special cases, we obtain new (or recover well-known) formulae for the Baum-Fulton-MacPherson Todd (or MacPherson-Chern) classes of toric varieties, as well as for the Thom-Milnor $L$-classes of simplicial projective toric varieties. 
We present two different perspectives for the computation of these characteristic classes of  toric varieties. First, we take advantage of the torus-orbit decomposition and  the motivic properties of the motivic Chern and resp. homology Hirzebruch classes to express the latter in terms of dualizing sheaves and resp. the (dual) Todd classes of closures of orbits. This method even applies to torus-invariant subspaces of a given toric variety. The obtained formula is then applied to  weighted lattice point counting in lattice polytopes and their subcomplexes, yielding generalized Pick-type formulae. Secondly, in the case of simplicial toric varieties, we compute our characteristic classes by using the Lefschetz-Riemann-Roch theorem of Edidin-Graham in the context of the geometric quotient description of such varieties. In this setting, we define mock Hirzebruch classes of  simplicial toric varieties (which specialize to the mock Chern,
mock Todd and mock $L$-classes of such varieties)
and investigate the difference  between the (actual) homology Hirzebruch class and the mock Hirzebruch class.  We show that this difference is localized on the singular locus, and we obtain a formula for it in which the contribution of each singular cone is identified explicitly. Finally, the two methods of computing characteristic classes are combined for proving several characteristic class formulae originally obtained by Cappell and Shaneson in the early $1990$s.
\end{abstract}

\maketitle

\tableofcontents

\section{Introduction}\label{intro}

Toric varieties admit several equivalent descriptions which arise naturally in many mathematical areas, including symplectic geometry, algebraic geometry, theoretical physics,  commutative algebra, etc. These perspectives, combined with intimate connections to pure and applied topics such as integer programming, representation theory, geometric modeling, algebraic topology, and 
combinatorics, lend toric varieties their importance, especially in view of their concreteness as examples.

\vspace{1.5mm}

A $d$-dimensional {\it toric variety} $X$ is an irreducible normal variety on which the complex $d$-torus acts with an open orbit, e.g., see \cite{CLS, Dan,Fu,O}. Toric varieties arise from combinatorial objects called {\it fans}, which are collections of cones in a lattice. Toric varieties are of interest both in their own right as algebraic varieties, and for their application to the theory of convex polytopes. For instance, Danilov \cite{Dan} used the Hirzebruch-Riemann-Roch theorem to establish a direct connection between the Todd classes of toric varieties and the problem of counting the number of lattice points in a convex polytope. Thus, the problem of finding explicit formulae for characteristic classes of toric varieties is of interest not only to topologists and algebraic geometers, but also to combinatorists, programmers, etc.

Let $X=X_{\Sigma}$ be the toric variety associated to a fan $\Sigma$ in a $d$-dimensional lattice $N$ of $\mathbb{R}^d$. We say that $X$ is {\it simplicial} if the fan is simplicial, i.e., each cone $\sigma$ of $\Sigma$ is spanned by linearly independent elements of the lattice $N$.  If, moreover, each cone $\sigma$ of $\Sigma$ is smooth (i.e., generated by a subset of a $\bZ$-basis for the lattice $N$), the corresponding toric variety $X_{\Sigma}$ is {\it smooth}.  Denote by $\Sigma(i)$ the set of $i$-dimensional cones of $\Sigma$. One-dimensional cones $\rho \in \Sig(1)$ are called {\it rays}. For each ray $\rho \in \Sig(1)$, we denote by $u_{\rho}$ the unique generator of the semigroup $\rho \cap N$. 
The $\{ u_{\rho} \}_{\rho \in \sig(1)}$ are the {\it generators} of $\sig$.
To each cone $\sigma \in \Sigma$ there corresponds a subvariety $V_{\sigma}$ of $X_{\Sigma}$, which is the closure of the orbit $O_{\sig}$ of $\sigma$ under the torus action. Each $V_{\sigma}$ is itself a toric variety of dimension $d - \dim(\sigma)$. 
The singular locus of $X_{\Sigma}$ is $\bigcup V_{\sigma}$, the union being taken over all singular cones $\sig \in \Sig_{\rm sing}$ in the fan $\Sigma$. Moreover, any toric variety 
is  Cohen-Macaulay and has only rational singularities. A simplicial toric variety has only (abelian) quotient singularities and is therefore a rational homology manifold
satisfying Poincar\'e duality over $\bQ$.

\vspace{1.5mm}

In this paper, we compute several characteristic classes of {\it (possibly singular) toric varieties}. More concretely, we compute  the {\it motivic Chern class} $mC_y(X_{\Sig})$, and the {\it homology Hirzebruch classes}, both the normalized 
$\widehat{T}_{y*}(X_{\Sig})$ and the un-normalized ${T}_{y*}(X_{\Sig})$, of such a toric variety $X_{\Sig}$. These characteristic classes depending on a parameter $y$ are defined in 
$G_0(X)[y]$ and resp. $H_*(X_{\Sig})\otimes \bQ[y]$, where $G_0(X)$ is the Grothendieck group of coherent sheaves, and for $H_*(-)$ one can use either the Chow group $CH_*$ 
or the even degree Borel-Moore homology group $H^{BM}_{2*}$.

By a result of Ishida \cite{I}, the motivic Chern class of a toric variety $X_{\Sig}$ is computed as \be\label{zero}
mC_y(X_{\Sig})=\sum_{p = 0}^d [\wti{\Omega}_{X_{\Sig}}^p] \cdot y^p \in G_0(X)[y], 
\ee
where the $\wti{\Omega}_{X_{\Sig}}^p$ denote the corresponding sheaves of Zariski differential forms. In fact, this formula holds even for a torus-invariant closed algebraic subset $X:=X_{\Sig'}\subseteq X_{\Sig}$ (i.e., a closed union of torus-orbits) corresponding to a star-closed subset $\Sig' \subseteq \Sig$, with the sheaves of differential $p$-forms  $\wti{\Omega}^p_{X}$ as introduced by Ishida \cite{I}.  Most of our results will apply to this more general context.

The un-normalized homology Hirzebruch class of such a torus-invariant closed algebraic subset $X:=X_{\Sig'}$ is defined as
\be\label{abis}{T}_{y*}(X):=td_* \circ mC_y(X)=\sum_{p = 0}^{d'} td_*([\wti{\Omega}_{X}^p]) \cdot y^p,\ee with $d'=\dim(X)$ and $td_*: G_0(-) \to H_*(-) \otimes \bQ$ the Todd class transformation of \cite{BFM, Fu2}. The normalized Hirzebruch class $\widehat{T}_{y*}(X)$ is obtained from ${T}_{y*}(X)$ by multiplying  by $(1+y)^{-k}$ on the degree $k$-part  ${T}_{y,k}(X) \in H_k(X)$, $k \geq 0$.

The degree of these characteristic classes for compact $X:=X_{\Sig'}$
recovers the $\chi_y$-genus of $X$ defined in terms of the Hodge filtration
of Deligne's mixed Hodge structure on $H_c^*(X,\bC)$:
\be\label{ione} 
\begin{split}
&\sum_{j,p}\,(-1)^j\dim_{\bC}\Gr_F^p H_c^j(X,\bC)\,(-y)^p =:\chi_y(X) \\
&=\sum_{p \geq 0} \chi(X,\wti{\Omega}_{X}^p) y^p =\int_{X}  T_{y*}(X)=\int_{X} \widehat{T}_{y*}(X).
\end{split}
\ee

The homology Hirzebruch classes, defined by Brasselet-Sch\"urmann-Yokura  in \cite{BSY},  have the virtue of unifying several other known characteristic class theories for singular varieties. For example, by specializing the parameter $y$ in the normalized Hirzebruch class $\widehat{T}_{y*}$ to the value $y=-1$, our results for Hirzebruch classes translate into formulae for the (rational) MacPherson-Chern classes $c_*(X)=\widehat{T}_{-1*}(X)$. By letting $y=0$, we get formulae for the Baum-Fulton-MacPherson Todd classes $td_*(X)=\widehat{T}_{0*}(X)$ of a torus-invariant closed algebraic subset $X$ of a toric variety, since by \cite{I}, such spaces only have Du Bois type singularities (i.e., $\cO_X \simeq \wti{\Omega}_X^0$). In particular, $\chi_0(X)=\chi(X,\cO_X)$ is the {\it arithmetic genus} in case $X$ is also compact. 
Moreover, for simplicial projective toric varieties or compact toric manifolds $X_{\Sig}$, we have for $y=1$ the following identification of the Thom-Milnor $L$-classes $$L_*(X_{\Sig})=\widehat{T}_{1*}(X_{\Sig}).$$ Conjecturally, the last equality should hold for a complete simplicial toric variety. The Thom-Milnor $L$-classes $L_*(X)\in H_{2*}(X,\bQ)$ are only defined for a compact rational homology manifold $X$, i.e., in our context for a complete simplicial toric variety.

\vspace{1.5mm}

We present two different perspectives for the computation of these characteristic classes of  toric varieties and their torus invariant closed algebraic subsets. 

\vspace{1.5mm}

\noindent{\bf Motivic Chern and Hirzebruch classes via orbit decomposition.}\\  
First, we take advantage of the torus-orbit decomposition of a toric variety together with  the motivic properties of these characteristic classes to express the motivic Chern and resp. homology Hirzebruch  classes in terms of dualizing sheaves and resp. the (dual) Todd classes of closures of orbits.  We prove the following result (see Theorem \ref{t1} and Theorem \ref{t2}):
\bt\label{ti1}
Let $X_{\Sig}$ be the toric variety defined by the fan $\Sigma$, and $X:=X_{\Sig'} \subseteq X_{\Sig}$ a  torus-invariant closed algebraic subset defined by a star-closed subset $\Sig' \subseteq \Sig$. For each cone $\sig \in \Sig$ with corresponding orbit $O_{\sig}$, denote by $k_{\sig} : V_{\sig} \hookrightarrow X$ the inclusion of the orbit closure.   Then the 
motivic Chern class $mC_*(X)$ is computed by:
\be\label{i0b} 
mC_y(X)=\sum_{\sig \in \Sig'}\, (1+y)^{\dim(O_{\sig})} \cdot (k_{\sig})_*[\omega_{V_{\sig}}],
\ee
and, similarly, the un-normalized homology Hirzebruch class ${T}_{y*}(X_{\Sig})$ is computed by the formula:
\be\label{i0} 
\begin{split}
T_{y*}(X)&=\sum_{\sig \in \Sig'}\, (1+y)^{\dim(O_{\sig})} \cdot (k_{\sig})_*td_*([\omega_{V_{\sig}}])\\
&=\sum_{\sig \in \Sig'}\, (-1-y)^{\dim(O_{\sig})} \cdot 
(k_{\sig})_*\left( td_*(V_{\sig})^{\vee}\right).
\end{split}
\ee
Here $\omega_{V_{\sig}}$ is the canonical sheaf of $V_{\sig}$ and the duality $(-)^{\vee}$
acts by multiplication by $(-1)^i$ on the $i$-th degree homology part $H_i(-)$.
The normalized Hirzebruch class $\widehat{T}_{y*}(X)$ is given by:
\be\label{i1}
\widehat{T}_{y*}(X)=\sum_{\sig\in \Sig', i}\, (-1-y)^{\dim(O_{\sig})-i}  \cdot (k_{\sig})_*td_i(V_{\sig}).
\ee
\et

These results should be seen as characteristic class versions of the well-known formula
\be\label{y}
\chi_y(X) =\sum _{\sig \in \Sig'}\, (-1-y)^{\dim(O_{\sig})}
\ee
for the $\chi_y$-genus of a toric variety, resp. its torus-invariant closed algebraic subsets, which one recovers for compact $X$
by taking the degree in (\ref{i0}) or (\ref{i1}) above, as 
$\int_{V_{\sig}} td_0(V_{\sig})=1$. For a far-reaching generalization of this $\chi_y$-genus formula to elliptic genera of complete toric manifolds, see \cite{BL}.

Letting, moreover, $y=-1,0,1$ in  (\ref{i0}) and (\ref{i1}), we obtain as special cases the following formulae:
\bc\label{ci1} (Chern, Todd and $L$-classes)\newline
Let $X_{\Sig}$ be the toric variety defined by the fan $\Sigma$, and $X:=X_{\Sig'} \subseteq X_{\Sig}$ a  torus-invariant closed algebraic subset defined by a star-closed subset $\Sig' \subseteq \Sig$. For each cone $\sig \in \Sig$ with corresponding orbit $O_{\sig}$, denote by $k_{\sig} : V_{\sig} \hookrightarrow X$ the inclusion of the orbit closure. 

\noindent $(a)$ (Ehler's formula) \ The (rational) MacPherson-Chern class $c_*(X)$ is computed by:
\be\label{i2}
c_*(X)=\sum_{\sig \in \Sig'}\,  (k_{\sig})_*td_{\dim(O_{\sig})}(V_{\sig})=\sum_{\sig \in \Sig'}\,  (k_{\sig})_*[V_{\sig}].
\ee
\noindent$(b)$ \ The Todd class $td_*(X)$ of $X$ is computed by:
\be\label{i3}
td_*(X)=\sum_{\sig \in \Sig'}\, (k_{\sig})_*td_*([\omega_{V_{\sig}}])
=\sum_{\sig \in \Sig'}\, (-1)^{\dim(O_{\sig})} \cdot 
(k_{\sig})_*\left( td_*(V_{\sig})^{\vee}\right).
\ee
\noindent$({c})$ \ The classes $T_{1*}(X)$
and  $\widehat{T}_{1*}(X)$ are computed by:
\be\label{i4un} 
T_{1*}(X)= \sum_{\sig \in \Sig'}\, 2^{\dim(O_{\sig})} \cdot (k_{\sig})_*td_*([\omega_{V_{\sig}}]) =
\sum_{\sig\in \Sig'}\, (-2)^{\dim(O_{\sig})} \cdot 
(k_{\sig})_*\left( td_*(V_{\sig})^{\vee}\right)
\ee
and
\be\label{i4nor} 
\widehat{T}_{1*}(X)=\sum_{\sig\in \Sig', i}\, (-2)^{\dim(O_{\sig})-i}  \cdot (k_{\sig})_*td_i(V_{\sig}).
\ee
Here $T_{1*}(X_{\Sig})=L^{AS}_*(X_{\Sig})$ is the Atiyah-Singer $L$-class for 
$X_{\Sig}$ a compact toric manifold, and $\widehat{T}_{1*}(X_{\Sig})=L_*(X_{\Sig})$
is the Thom-Milnor $L$-class for $X_{\Sig}$ a compact toric manifold
or a simplicial projective toric variety.
\ec
For the last identification of $\widehat{T}_{1*}$-classes with the Thom-Milnor $L$-classes we make use of a recent result of \cite{GS} which asserts that simplicial projective toric varieties can be realized as global quotients of a (projective) manifold by a finite algebraic group action,
so that we can apply a corresponding result from \cite{CMSS}. Conjecturally, this identification should hold for a complete simplicial toric variety.
Taking the degree in (\ref{i4nor}), we get the following description of the signature
$sign(X_{\Sig})$ of a compact toric manifold
or a simplicial projective toric variety $X_{\Sig}$
(compare also with \cite{LR}, and in the smooth case also with \cite{O}[Thm.3.12]):
\be
sign(X_{\Sig})=\chi_1(X_{\Sig})=\sum _{\sig \in \Sig}\, (-2)^{\dim(O_{\sig})}.
\ee
Also note that Ehler's formula
for the MacPherson-Chern class $c_*(X_{\Sig})$ of a toric variety $X_{\Sig}$ (even with integer coefficients) is due to
\cite{Al, BBF}. Theorem \ref{ti1} should be seen as a counterpart of Ehler's formula for the motivic Chern and Hirzebruch classes. In fact, our proof of Theorem \ref{ti1} follows closely the strategy of Aluffi \cite{Al} for his proof of  Ehler's formula. Formula
(\ref{i3}) for the Todd class $td_*(X)$ of a torus-invariant closed algebraic subset of a toric variety $X_{\Sig}$ is (to our knowledge) completely new, and is obtained here from the motivic nature of the Hirzebruch classes via the identification $td_*(X)=T_{0*}(X)$ (which uses the fact that such a space has only Du Bois singularities, cf. \cite{I}).

\vspace{1.5mm}

Moreover, in the case of a toric variety $X_{\Sig}$ (i.e., $\Sig'=\Sig$), one gets from (\ref{i3}) and (\ref{i4un}) by using a ``duality argument''  also the 
following formulae:
\be\label{i3dual}
td_*([\omega_{X_{\Sig}}])=
\sum_{\sig \in \Sig}\, (-1)^{{\rm codim}(O_{\sig})} \cdot 
(k_{\sig})_* td_*(V_{\sig})
\ee
and in the simplicial context (see Proposition \ref{4Dual}):
\be\label{i4dual} 
T_{1*}(X_{\Sig})= (-1)^{\dim(X_{\Sig})}\cdot \left(T_{1*}(X_{\Sig})\right)^{\vee}= 
\sum_{\sig\in \Sig}\, (-1)^{\dim(X_{\Sig})}\cdot (-2)^{\dim(O_{\sig})} \cdot 
(k_{\sig})_* td_*(V_{\sig}).
\ee

The above formulae calculate the homology Hirzebruch classes (and, in particular, the Chern,  Todd and $L$-classes, respectively) of singular toric varieties in terms of the Todd classes of closures of torus orbits. Therefore, known results about Todd classes of toric varietes (e.g., see  \cite{BP, BeV, BV, BZ, GP, Mor, P, P2, P3, PT}) also apply to these other characteristic classes: in particular, Todd class formulae satisfying the so-called {\it Danilov condition}, as well as the algorithmic computations from \cite{BP, P3}. 
As it  will be discussed later on, in the simplicial context  these  Todd classes can also be computed by using the Lefschetz-Riemann-Roch theorem of \cite{EG} for geometric quotients.

\vspace{1.5mm}

\noindent{\bf Weighted lattice point counting.}\\
Let $M$ be a lattice and $P \subset M_{\bR} \cong \bR^d$ be a full-dimensional lattice polytope with associated projective toric variety $X_P$ and ample Cartier divisor $D_P$.
By the classical work of Danilov \cite{Dan},
the (dual) Todd classes of $X_P$ can be used for counting the number of lattice points in (the interior of) a lattice polytope $P$:
\begin{equation*}
\begin{split}
| P \cap  M | &=\int_{X_P} ch(\cO_{X_P}(D_P)) \cap td_*(X_P) =\sum_{k\geq 0} \frac{1}{k!} \int_{X_P} [D_P]^k \cap td_k(X_P), \\
| \Int({P}) \cap M | &=\int_{X_P} ch(\cO_{X_P}(D_P)) \cap td_*([\omega_{X_P}])=\sum_{k\geq 0} \frac{(-1)^d}{k!} \int_{X_P} [-D_P]^k \cap td_k(X_P),
\end{split}
\end{equation*}
as well as the coefficients of the {\it Ehrhart polynomial of $P$} counting the number of lattice points in the dilated polytope $\ell P:=\{\ell \cdot u \ | \ u \in P \}$ for a positive integer $\ell$ (see Section \ref{countsec} for details):
$${\rm Ehr}_P(\ell):=|\ell P \cap M|=\sum_{k=0}^d a_k \ell^k,$$ with $$a_k= \frac{1}{k!}  \int_{X_P}  [D_P]^k \cap td_k(X_P).$$

In the special case of a lattice polygon $P \subset \bR^2$, the Todd class of the corresponding surface $X_P$ is given by the well-known formula (e.g., see \cite{Fu}[p.130]) 
\be\label{pol}
td_*(X_P)=[X_P]+\frac{1}{2}\sum_{\rho \in \Sig_P(1)} [V_{\rho}] + [{\rm pt}],
\ee
where $\Sig_P$ denotes the fan associated to $P$. So by evaluating the right-hand side of the above lattice point counting formula for $P$, one gets the classical {\it Pick's formula}:
\be\label{Pick}
| P \cap  M | ={\rm Area}({P}) + \frac{1}{2} | \partial P \cap M | +1,
\ee
where ${\rm Area}({P})$ is the usual Euclidian area of $P$.

As a direct application of Theorem \ref{ti1}, we show that the un-normalized homology Hirzebruch classes are useful for lattice point counting with certain weights, reflecting the face decomposition
$$P = \bigcup_{Q \preceq P}\: \Relint(Q).$$
In more detail, we prove the following {\it generalized Pick-type formulae} (see Theorem \ref{wc}, and compare also with \cite{Ag} and \cite{M}):
\bt\label{ti2} Let $M$ be a lattice of rank $d$ and $P \subset M_{\bR} \cong \bR^d$ be a full-dimensional lattice polytope with associated projective toric variety $X_P$ and ample Cartier divisor $D_P$. In addition, let $X:=X_{P'}$ be a torus-invariant closed algebraic subset of $X_P$ corresponding to a polytopal subcomplex $P' \subseteq P$ (i.e., a closed union of faces of $P$). Then the following formula holds:
\be\label{i5}
\sum_{Q \preceq P'} (1+y)^{\dim(Q)} \cdot | \Relint(Q) \cap M |=\int_{X} ch(\cO_{X_P}(D_P)|_X) \cap T_{y*}(X),
\ee
where the summation on the left is over the faces $Q$ of $P'$, $|-|$ denotes the cardinality of sets and $\Relint(Q)$ the relative interior of a face $Q$.
\et
\bc
The Ehrhart polynomial of the polytopal subcomplex $P' \subseteq P$ as above, counting the number of lattice points in the dilated complex $\ell P':=\{\ell \cdot u \ | \ u \in P' \}$ for a positive integer $\ell$, is computed by:
$${\rm Ehr}_{P'}(\ell):=|\ell P' \cap M|=\sum_{k=0}^{d'} a_k \ell^k,$$ with $$a_k= \frac{1}{k!}  \int_{X}  [D_P|_X]^k \cap td_k(X),$$
and $X:=X_{P'}$ the corresponding  $d'$-dimensional torus-invariant closed algebraic subset of $X_P$.
\ec

Note that ${\rm Ehr}_{P'}(0)=a_0=\int_X td_*(X)$ can be computed on the combinatorial side by formula (\ref{y}) as the Euler characteristic of $P'$, while on the algebraic geometric side it is given by the arithmetic genus of $X:=X_{P'}$: 
\be\label{Er} \chi({P'}):=\sum_{Q \preceq P'} (-1)^{\dim Q}=\chi_0(X)=\int_X td_*(X)=\chi(X,\cO_X).\ee

Moreover, by using the duality formula (\ref{i4dual}), we also have (see formula (\ref{5L})):
\bt
If $X_P$ is the simplicial projective toric variety associated to a simple polytope $P$, with corresponding  ample Cartier divisor $D_P$, then the following formula holds:
\be\label{i5L}
\sum_{Q \preceq P} \left(-\frac{1}{2}\right)^{{\rm codim}(Q)} \cdot |Q \cap M |=\int_{X_P} ch(\cO_{X_P}(D_P)) \cap \left( \frac{1}{2}\right)^{\dim(X_P)}\cdot T_{1*}(X_P) .
\ee
\et

Note that if $X_P$ is smooth (i.e., $P$ is a Delzant lattice polytope), the un-normalized
Hirzebruch class $T_{1*}(X_P)$ corresponds under Poincar\'e duality to the characteristic
class of the tangent bundle defined in terms of Chern roots by the un-normalized power series
$\alpha/\tanh(\frac{\alpha}{2})$, so that $\left( \frac{1}{2}\right)^{\dim(X_P)}\cdot T_{1*}(X_P)$
corresponds to the characteristic
class of the tangent bundle defined in terms of Chern roots  by the normalized power series
$\frac{\alpha}{2}/\tanh(\frac{\alpha}{2})$. 

Therefore (\ref{i5L}) extends the result of 
\cite{Fe}[Corollary 5] for Delzant lattice polytopes to the more general case of simple lattice polytopes. Our approach is very different from the techniques used in \cite{Fe} for projective toric manifolds, which are based on virtual submanifolds Poincar\'e dual to cobordism characteristic classes of the projective toric manifold.

\bex {\it Parametrized Pick's formula for lattice polygones} \\
In the special case of a lattice polygon $P \subset \bR^2$ (which  is always simple, by definition), the un-normalized homology Hirzebruch class  $T_{y*}(X_P)$ of the corresponding projective surface $X_P$ is given by (see Remark \ref{rem_mock})
\be\label{pol}
T_{y*}(X_P)=(1+y)^2 \cdot [X_P]+\frac{1-y^2}{2}\sum_{\rho \in \Sig_P(1)} [V_{\rho}] + \chi_y(X_P) \cdot [{\rm pt}],
\ee
where $\Sig_P$ denotes as before the fan associated to $P$. So, by evaluating the right-hand side of the above weighted lattice point counting formula (\ref{i5}), we obtain the following parametrized version of the classical Pick's formula (\ref{Pick}):
\be\label{yPick}
\begin{split}
(1+y)^2 \cdot | {\rm Int}({P}) \cap  M | + (1+y) \cdot \sum_{F} |  \Relint(F) \cap M | + v \\ =(1+y)^2 \cdot {\rm Area}({P}) + \frac{1-y^2}{2} | \partial P \cap M | +\chi_y({P}),
\end{split}
\ee
where $F$ runs over the facets (i.e., boundary segments) of $P$,  $v$ denotes the number of vertices of $P$, and $$\chi_y({P}):=(1+y)^2 - (1+y) \cdot | \{ F \} | + v=\chi_y(X_P).$$
By specializing (\ref{yPick}) to $y=0$, we recover Pick's formula (\ref{Pick}), since $\chi_0({P})=1$.
Similarly, by letting $y=1$ and dividing by $4$, formula (\ref{yPick}) translates into:
\begin{equation*}
\begin{split}
| {\rm Int}({P}) \cap  M | + \frac{1}{2} \cdot \sum_{F} |  \Relint(F) \cap M | + \frac{v}{4}  = {\rm Area}({P}) + \frac{4- v}{4},
\end{split}
\end{equation*}
with $4- v=\chi_1({P})=sign(X_P)$. 
This is easily seen to be equivalent to Pick's formula (\ref{Pick}), compare also with \cite{Fe}[Introduction].
Moreover, by using formula (\ref{i5L})  we get the following identity (see also \cite{Fe}[eqn.(7)]):
\begin{equation*}
\begin{split}
| P \cap  M | - \frac{1}{2} \cdot \sum_{F} |  F \cap M | + \frac{v}{4}  = {\rm Area}({P}) + \frac{4- v}{4},
\end{split}
\end{equation*}
which is yet another re-writing of Pick's formula. Finally, it can be easily seen that the parametrized Pick formula (\ref{yPick}) is itself equivalent to the classical Pick formula (\ref{Pick}).
\eex

\vspace{1.5mm}

\noindent{\bf Generalized toric Hirzebruch-Riemann-Roch.}\\
By using a well-known ``module property'' of the Todd class transformation $td_*$, one obtains the following:
\bt\label{i-gHRR}(generalized toric Hirzebruch-Riemann-Roch)\newline
Let $X_{\Sig}$ be the toric variety defined by the fan $\Sigma$, and $X:=X_{\Sig'} \subseteq X_{\Sig}$ be a torus-invariant closed algebraic subset of dimension $d'$ defined by a star-closed subset $\Sig' \subseteq \Sig$. Let $D$ be a fixed Cartier divisor on $X$. Then the (generalized) {\it Hirzebruch polynomial of $D$}:
\be
\chi_y(X,\cO_X(D)):=\sum _{p=0}^{d'} \chi(X,\wti{\om}^p_X\otimes \cO_X(D)) \cdot y^p,
\ee
with $\wti{\om}^p_X$ the corresponding Ishida sheaf of $p$-forms, 
 is computed by the formula:
\be\label{i-Hp}
\chi_y(X,\cO_X(D))=\int_{X} ch(\cO_X(D)) \cap T_{y*}(X).
\ee
\et
By combining the above formula with Theorem \ref{ti2}, we get the following combinatorial description
(compare with \cite{M} for another approach to this result in the special case $\Sig'=\Sig$):

\bc\label{chiy} Let $M$ be a lattice of rank $d$ and $P \subset M_{\bR} \cong \bR^d$ be a full-dimensional lattice polytope with associated projective simplicial toric variety $X=X_P$ and ample Cartier divisor $D_P$. In addition, let $X:=X_{P'}$ be a torus-invariant closed algebraic subset of $X_P$ corresponding to a polytopal subcomplex $P' \subseteq P$ (i.e., a closed union of faces of $P$). Then the following formula holds:
\be
\chi_y(X,\cO_X(D_P|_X)) = \sum_{Q \preceq P'} (1+y)^{\dim(Q)} \cdot | \Relint(Q) \cap M |.
\ee
\ec

\vspace{1.5mm}

\noindent{\bf Hirzebruch classes of simplicial toric varieties via Lefschetz-Riemann-Roch.}\\
It is known that simplicial toric varieties have an intersection theory, provided that $\bQ$-coefficients are used. The generators of the rational cohomology (or Chow) ring are the classes $[V_{\rho}]$ defined by the $\bQ$-Cartier divisors corresponding to the rays $\rho \in \Sig(1)$ of the fan $\Sig$. So it natural to try to express all characteristic classes in the simplicial context in terms of these generators. This brings us to our second computational method, which employs the Lefschetz-Riemann-Roch theorem  for the geometric quotient realization $X_{\Sig}=W/G$ of such simplicial toric varieties (without a torus factor). 
In fact, here one can also use other such quotient presentations of $X_{\Sig}$ in terms of a {\it stacky fan} 
in the sense of Borison-Chen-Smith \cite{BCS}.
In the case of the Cox construction, this approach was already used by Edidin-Graham \cite{EG} for the computation of Todd classes of simplicial toric varieties, extending a corresponding result of  Brion-Vergne \cite{BV}
for complete simplicial toric varieties (compare also with \cite{G} for a related approach using Kawasaki's orbifold Riemann-Roch theorem, with applications to Euler-MacLaurin formulae).

\vspace{1.5mm}

In this paper, we enhance the calculation of  Edidin-Graham \cite{EG} to the homology Hirzebruch classes,
based on the description (\ref{abis}) of the un-normalized Hirzebruch classes in terms of 
the sheaves of Zariski $p$-forms (see Theorem \ref{Tc}):
 
\bt\label{Tc-intro}  
The normalized Hirzebruch class of a simplicial toric variety $X=X_{\Sig}$ 
 is given by:
\be\label{nHirzc-intro}
\widehat{T}_{y*}(X)= \left( \sum_{g \in G_{\Sig}} \: \prod_{\rho \in \Sig(1)} \frac{[V_{\rho}] \cdot 
\big( 1+y  \cdot a_{\rho}(g)  \cdot e^{-[V_{\rho}](1+y)}\big)}{1-a_{\rho}(g) \cdot e^{-[V_{\rho}](1+y)}} \right) \cap [X]\:.
\ee
\et
For more details, especially for the definition of $G_{\Sig}$ and $a_{\rho}(g)$, we refer to 
Theorem \ref{Tc}. Here we only point out that $G_{\Sig} \subset G$ is the set of group elements having fixed-points in the Cox quotient realization $X_{\Sig}=W/G$. 

\bex {\it Weighted projective spaces} \\
Consider the example of a $d$-dimensional weighted projective space $X_{\Sig}=\bP(q_0,\dots,q_d):=(\bC^{d+1} \setminus \{0\})/G$, with $G \simeq \bC^*$ acting with weights $(q_0,\dots,q_d)$. By \cite{RT}[Thm.1.26], we can assume that the weights are {\it reduced}, i.e., $1=\gcd(q_0,\dots , q_{j-1}, q_{j+1},\dots q_d)$, for all $j=0,\dots,d$, so that this quotient realization is the Cox construction. The normalized Hirzebruch class is then computed by (see Example \ref{fake}):
\be\label{wph0}
\widehat{T}_{y*}(X_{\Sig})= \left( \sum_{\lambda \in G_{\Sig}} \: \prod_{j=0}^d \frac{q_j T \cdot 
\big( 1+y  \cdot \lambda^{q_j}  \cdot e^{-q_j (1+y) T}\big)}{1-\lambda^{q_j}  \cdot e^{-q_j (1+y) T}} \right) \cap [X_{\Sig}]\:,
\ee
with $$G_{\Sig}:=\bigcup_{j=0}^d \{ \lambda\in \bC^* \ | \ \lambda^{q_j}=1 \}$$ 
and $T=\frac{[V_{j}]}{q_j}\in H^2(X_{\Sig};\bQ)$, $j=0,\dots, d$, the distinguished generator of the rational cohomology ring $H^*(X_{\Sig};\bQ)\simeq \bQ[T]/(T^{d+1}).$ 
The un-normalized Hirzebruch class of $X_{\Sig}$ is obtained from the above formula (\ref{wph0}) by substituting $T$ for $(1+y)T$ and $(1+y)^d \cdot [X] $ in place of $[X]$. The resulting formula recovers Moonen's calculation in \cite{M}[p.176] (see also \cite{CMSS}[Ex.1.7]).
\eex

By specializing to $y=0$, one  recovers from (\ref{nHirzc-intro}) exactly the result of Edidin-Graham \cite{EG} for the Todd class $td_*(X)=\widehat{T}_{0*}(X)$ of $X$.
Moreover, in the case of a smooth toric variety one has $G_{\Sig}=\{id\}$, with all $a_{\rho}(id)=1$.
In particular, we also obtain formulae for the Hirzebruch classes of smooth toric varieties, in which the generators $[V_{\rho}]$ ($\rho \in \Sig(1)$) behave like the Chern roots of the tangent bundle of $X_{\Sig}$ (except that we may have $|\Sig(1)|>\dim(X_{\Sig})$):

\bc\label{sm-intro}
The normalized Hirzebruch class of a smooth toric variety $X=X_{\Sig}$ is given by:
\be\label{Hirzsm-intro}
\widehat{T}_{y*}(X)= \left( \prod_{\rho \in \Sig(1)} \frac{[V_{\rho}] \cdot \big( 1+y \cdot e^{-[V_{\rho}](1+y)}\big)}{1-e^{-[V_{\rho}](1+y)}} \right)
\cap [X] \:.
\ee
\ec

\vspace{1.5mm}

\noindent{\bf Mock Hirzebruch classes.}\\
Following the terminology from \cite{CS,P}, we define the (normalized)  {\it mock Hirzebruch class} $\widehat{T}^{(m)}_{y*}(X)$ of a {simplicial} toric variety $X=X_{\Sigma}$ to be the class computed using the formula (\ref{Hirzsm-intro}). Similarly, by taking the appropriate polynomials in the $[V_{\rho}]$'s, one can define mock Todd classes (see \cite{P}), 
mock Chern classes and mock $L$-classes, all of which are specializations 
of the mock Hirzebruch class obtained by evaluating the parameter $y$ at $-1$, $0$, and $1$, respectively. 

In Section \ref{M},  we investigate the difference $\widehat{T}_{y*}(X)-\widehat{T}^{(m)}_{y*}(X) $ between the actual (normalized) Hirzebruch class and the mock Hirzebruch class 
of a simplicial toric variety $X=X_{\Sig}$. We show that this difference is localized on the singular locus, and we obtain a formula for it in which the contribution of each singular 
cone $\sig \in \Sig_{\rm sing}$  is identified explicitly. We prove the following (see Theorem \ref{Tcb} and Remark \ref{rem_mock}):
\bt\label{ti3} The normalized homology Hirzebruch class of a simplicial toric variety $X=X_{\Sig}$ is computed by the formula:
\be\label{i6}
\widehat{T}_{y*}(X)=\widehat{T}^{(m)}_{y*}(X) + \sum_{\sig \in \Sig_{\rm sing}}  \cA_y({\sig}) \cdot \left( (k_{\sig})_*\widehat{T}^{(m)}_{y*}(V_{\sig}) \right) \:,
\ee
with $k_{\sig}:V_{\sig} \hookrightarrow X$ denoting the orbit closure inclusion,
and 
\be\label{A-intro} \cA_y({\sig}):= \frac{1}{\mult(\sig)} \cdot \sum_{g \in G_{\sig}^{\circ}}\prod_{\rho \in \sig(1)} 
\frac{ 1+y  \cdot a_{\rho}(g)  \cdot e^{-[V_{\rho}](1+y)}}{1-a_{\rho}(g) \cdot e^{-[V_{\rho}](1+y)}}.\ee
Here, the numbers $a_{\rho}(g)$ are roots of unity of order $\mult(\sig)$, different from $1$.
In particular, since $X$ is smooth in codimension one, we get:
\be\label{i6b}
\widehat{T}_{y*}(X)=[X] + \frac{1-y}{2} \sum_{\rho \in \Sig(1)} [V_{\rho}] + {\rm lower \ order\ homological \ degree \ terms.}
\ee
\et
The corresponding mock Hirzebruch classes $\widehat{T}^{(m)}_{y*}(V_{\sig})$ of orbit closures in (\ref{i6})  can be regarded as {\it tangential data} for the fixed point sets of the action in the geometric quotient description of $X_{\Sig}$, whereas the coefficient $\cA_y(\sig)$ encodes the {\it normal data} information. 

It should be noted that in the above coefficient $\cA_y(\sig)$ only the cohomology classes $[V_{\rho}]$ depend on the fan $\Sig$, but the multiplicity $\mult(\sig)$ and the character $a_{\rho}:G_{\sig} \to \bC^*$ depend only on rational simplicial cone $\sig$ (and not on the fan $\Sig$ nor the group $G$ from the Cox quotient construction). In fact, these can be given directly as follows.
For a $k$-dimensional rational simplicial cone $\sig$ generated by the rays $\rho_1,\dots,\rho_k$ one defines
the finite abelian group $G_{\sig}:=N_{\sig}/(u_1,\dots,u_k)$  as the quotient of the sublattice $N_{\sig}$ of $N$ spanned by the points in $\sig\cap N$ modulo
 the sublattice $(u_1,\dots,u_k)$ generated by the ray generators $u_j$ of  $\rho_j$.
Then $|G_{\sig}|=\mult(\sig)$ is just the multiplicity of $\sig$, with $\mult(\sig)=1$ exactly in case of a smooth cone.
 Let $m_i\in M_{\sig}$ for $1\leq i\leq k$ be the unique primitive elements in the dual lattice $M_{\sig}$ of $N_{\sig}$ satisfying $ \langle m_i,u_j \rangle = 0$ for $i\neq j$ and
 $ \langle m_i,u_i \rangle > 0$, so that the dual lattice $M_{\sig}'$ of $(u_1,\dots,u_k)$ is generated by the elements $\frac{m_j}{\langle m_j,u_j \rangle}.$
  Then $a_{\rho_j}(g)$ for $g=n+(u_1,\dots, u_k)\in G_{\sig}$ is defined as 
 $$a_{\rho_j}(g):=\exp\left(2\pi i\cdot \gamma_{\rho_j}(g)
 \right), \quad \text{with} \quad \gamma_{\rho_j}(g):= 
 \frac{\langle m_j,n \rangle}{\langle m_j,u_j \rangle}\:.  $$
Moreover, 
 $$G_{\sig}^{\circ}:=\{g\in G_{\sig}|\: a_{\rho_j}(g)\neq 1 \:\text{for all $1\leq j\leq k$}\}.$$
 Finally, using the parallelotopes  $$P_{\sig}:=\{ \sum_{i=1}^k \lambda_i u_i | \ 0 \leq \lambda_i <1 \} \quad {\rm and} \quad 
 P^{\circ}_{\sig}:=\{ \sum_{i=1}^k \lambda_i u_i | \ 0 < \lambda_i <1 \},$$ 
we also have 
  $G_{\sig}^{\circ} \simeq N \cap  P^{\circ}_{\sig}$ via the natural bijection  $G_{\sig} \simeq N \cap  P_{\sig}$.
 
 \vspace{1.5mm}
 
 Note that for $y=-1$ one gets by $\cA_{-1}(\sig)=\frac{|G_{\sig}^{\circ}|}{\mult(\sig)}$
 an easy relation between Ehler's formula (\ref{i2}) for the MacPherson Chern class and the mock
 Chern class of a toric variety which is given by
 $$c_*^{(m)}(X_{\Sig}):= \left( \prod_{\rho \in \Sig(1)} (1+ [V_{\rho}])  \right) \cap [X_{\Sig}] 
= \sum_{\sig\in \Sig}\: \frac{1}{\mult(\sig)}\cdot [V(\sig)].$$
 
 Similarly, the top-dimensional contribution of a singular cone $\sig \in \Sig_{\rm sing}$ in  formula
(\ref{i6}) is computed by
\be\label{topcont}
\frac{1}{\mult(\sig)} \cdot \sum_{g \in G_{\sig}^{\circ}}  \left( \prod_{\rho \in \sig(1)} 
\frac{ 1+y  \cdot a_{\rho}(g)}{1-a_{\rho}(g)}\right) \cdot [V_{\sig}],
\ee
 which for an isolated singularity $V_{\sig}=\{pt\}$ is already the full contribution.
 For a singular cone of smallest codimension, we have that  $G_{\sig}^{\circ}=G_{\sig}\backslash \{id\}$ and $G_{\sig}$ is cyclic of order $\mult(\sig)$, so that  $a_{\rho}(g)$  ($g \in G_{\sig}^{\circ}$) are  primitive roots of unity.
 Specializing further to $y=0$ (resp. $y=1$) we recover the correction factors
 for Todd (and resp. $L$-) class formulae of toric varieties appearing in  \cite{BP, BV, GP, P, P2}) (resp. \cite{Za} in the case of weighted projective spaces) also in close relation to (generalized) Dedekind (resp. cotangent) sums. 
 
Moreover,  for  $y=1$, we get for  a projective simplicial toric variety $X_{\Sig}$  the correction terms of the singular cones for the difference between the Thom-Milnor  and resp. the mock
$L$-classes of $X_{\Sig}$
(where the mock $L$-classes are used implicitly in \cite{LR}). In particular, the equality
$$sign(X_{\Sig}) = \int_{X_{\Sig}} L^{(m)}_*(X_{\Sig})$$
used in \cite{LR} and \cite{Fe}[Ex.1]  holds in general only for a smooth projective toric variety. 
\bex\label{ex-intro} {\it Weighted projective spaces, revisited.}\\
Consider the example of a $d$-dimensional weighted projective space $X_{\Sig}=\bP(1, \dots, 1,m)$ with weights $(1, \dots, 1,m)$ and $m >1$. This is a projective simplicial toric variety which has exactly one isolated singular point for $d \geq 2$. Formula (\ref{i6}) becomes in this case (see Example \ref{exw}):
\be
\widehat{T}_{y*}(X_{\Sig})=\widehat{T}^{(m)}_{y*}(X_{\Sig}) + \frac{1}{m} \sum_{ \lambda^m=1,\\  \lambda \neq 1}  \left( \frac{1+\lambda y}{1-\lambda} \right)^d \cdot [pt]\:.
\ee
In particular, for $y=1$, we get:
$$sign( X_{\Sig})=\int_{X_{\Sig}}L^{(m)}_*(X_{\Sig}) + \frac{1}{m} \sum_{ \lambda^m=1,\\  \lambda \neq 1}  \left( \frac{1+\lambda}{1-\lambda} \right)^d \:,$$
and it can be easily seen that the correction factor does not vanish in general, e.g., for $d=2$, the last sum is given by $$\sum_{ \lambda^m=1,\\  \lambda \neq 1}  \left( \frac{1+\lambda}{1-\lambda} \right)^2=-\frac{1}{3}(m-1)(m-2),$$ see \cite{Ge}[p.7]. Similarly, for $y=0$, we obtain for the arithmetic genus 
$$\chi_0( X_{\Sig})=\int_{X_{\Sig}}td^{(m)}_*(X_{\Sig}) + \frac{1}{m} \sum_{ \lambda^m=1,\\  \lambda \neq 1}  \left( \frac{1}{1-\lambda} \right)^d \:,$$
where the last sum is given for $d=2$ by:
$$\sum_{ \lambda^m=1,\\  \lambda \neq 1}  \left( \frac{1}{1-\lambda} \right)^2=-\frac{(m-1)(m-5)}{12},$$ see \cite{Ge}[p.4].
\eex


\vspace{1.5mm}

\noindent{\bf The characteristic class formulae of Cappell-Shaneson.}\\
In Section \ref{CS90}, we combine the two computational perspectives described above for proving several formulae originally obtained by Cappell and Shaneson in the early 1990s,
see \cite{CS,S}.

The following formula (in terms of $L$-classes) appeared first in Shaneson's $1994$ ICM proceeding paper, see \cite{S}[(5.3)].
\bt\label{CS1-intro} 
Let the {\it $T$-class} of a $d$-dimensional simplicial toric variety $X$  be defined as
\be
T_*(X):=\sum_{k=0}^d 2^{d-k} \cdot td_k(X).
\ee
For any cone $\sig \in \Sig$, let $k_{\sig}:V_{\sig} \hookrightarrow X$ be the inclusion of the corresponding orbit closure. Then the following holds:
\be\label{cs1-intro}
T_*(X)=\sum_{\sig \in \Sig} (k_{\sig})_*\widehat{T}_{1*}(V_{\sig}).
\ee
\et
Recall here that the identification $L_*(X_{\Sig})=\widehat{T}_{1*}(X_{\Sig})$
holds for a simplicial projective toric variety or a compact toric manifold. Conjecturally, this equality should hold for any complete simplicial toric variety. 

While Cappell-Shaneson's method of proof of the above result uses mapping formulae for Todd and $L$-classes in the context of resolutions of singularities, our approach relies on Theorem \ref{ti1}, i.e., on the additivity properties of the motivic Hirzebruch classes.

\vspace{1.5mm}

The above theorem suggests the following definition of mock $T$-classes (cf. \cite{CS,S}) in terms of
mock $L$-classes $L_*^{(m)}(X)=\widehat{T}^{(m)}_{1*}(X)$ of a  simplicial toric variety 
$X=X_{\Sig}$ (and similarly for the orbit closures $V_\sig$):
\be\label{mockt-intro}
T_*^{(m)}(X):=\sum_{\sig \in \Sig} (k_{\sig})_*L_*^{(m)}(V_{\sig}).
\ee

Then, by using Theorem \ref{CS1-intro}, we have as in \cite{CS}[Theorem 4] the following expression for the $T$-classes (which are, in fact, suitable renormalized Todd classes) in terms of  mock $T$- and resp. mock $L$-classes,
with $ \cA_y({\sig})$ as in (\ref{A-intro}):
\bt\label{CS2-intro} Let $X=X_{\Sig}$ be a simplicial toric variety. Then in the above notations we have that:
\be\label{cs2-intro}
T_*(X) = T_*^{(m)}(X)+\sum_{\sig \in \Sig_{\rm sing}}\:\cA_{1}({\sig})
\cdot (k_{\sig})_*T_*^{(m)}(V_{\sig}) \:.
\ee
\et
While Cappell-Shaneson's approach to the above result uses induction and Atiyah-Singer (resp. Hirzebruch-Zagier) type results for $L$-classes in the {\it local} orbifold description of toric varieties, we use the specialization of Theorem \ref{ti3} to the value $y=1$, which is based on Edidin-Graham's Lefschetz-Riemann-Roch theorem in the context of the Cox {\it global} geometric quotient construction.

Finally, the following renormalization of the result of Theorem \ref{CS2-intro} in terms of coefficients
\be\label{alpha-intro} \begin{split} \alpha(\sig)&:= \frac{1}{\mult(\sig)} \cdot \sum_{g \in G_{\sig}^{\circ}}\prod_{\rho \in \sig(1)} 
\frac{ 1+ a_{\rho}(g)  \cdot e^{-[V_{\rho}]}}{1-a_{\rho}(g) \cdot e^{-[V_{\rho}]}}\\
&= \frac{1}{\mult(\sig)} \cdot \sum_{g \in G_{\sig}^{\circ}}\prod_{\rho \in \sig(1)} 
\coth \left( \pi i \cdot \gamma_{\rho}(g)+\frac{1}{2}[V_{\rho}] \right)
\end{split}
\ee
fits better with the corresponding Euler-MacLaurin formulae (see \cite{CS2}[Thm.2], \cite{S}[Sect.6]), with
$\alpha(\{0\}):=1$ and $\alpha(\sig):=0$ for any other smooth cone $\sig\in \Sig$.

\bc Let $X=X_{\Sig}$ be a simplicial toric variety. Then we have:
\be\label{todd-EM-intro}
td_*(X)= \sum_{\sig\in \Sig}\: \alpha(\sig)\cdot \left(
\sum_{\{\tau | \sig \preceq \tau\}}\: \mult(\tau)\prod_{\rho \in \tau(1)} \frac{1}{2}[V_{\rho}]
\prod_{\rho \notin \tau(1)} \frac{\frac{1}{2}[V_{\rho}]}{\tanh(\frac{1}{2}[V_{\rho}])}
\right) \cap [X].
\ee
\ec

However,  in this paper we only focus on geometric results for characteristic classes of
(simplicial) toric varieties. The relation with the corresponding Euler-MacLaurin formulae will be discussed elsewhere.

\begin{ack} 
We thank Sylvain Cappell and Julius Shaneson for encouraging us to pursue this project, and for explaining to us 
their approach towards computing Todd and $L$-classes of toric varieties with applications to Euler-Maclaurin formulae.  L. Maxim is partially supported by NSF-1005338 and by a research fellowship from the Max-Planck-Institut f\"ur Mathematik, Bonn.
J. Sch\"urmann is supported by the SFB 878 ``groups, geometry and actions".

\end{ack}


\section{Preliminaries}

\subsection{Motivic Chern and Homology Hirzebruch classes of singular varieties}\label{motHirs}
In this section, we recall the construction and main properties of the motivic Chern and Hirzebruch classes of {\it singular} complex algebraic varieties, as developed in \cite{BSY} (but see also \cite{Sch}).

\subsubsection{Motivic Chern and Hirzebruch classes}\label{motHir}
For a given complex algebraic variety $X$, let $K_{0}(var/X)$ be the  relative Grothendieck group of complex algebraic varieties over $X$ as introduced and studied by Looijenga \cite{Lo}
and Bittner \cite{Bi} in relation to motivic integration.
Here a complex variety is a separated scheme of finite type over $spec(\bC)$.
$K_{0}(var/X)$ is the quotient of the free abelian group of isomorphism
classes of algebraic morphisms $Y\to X$ to $X$,  modulo the ``additivity''
relation generated by  
$$
[Y\to X] = [Z\to Y \to X] + [Y\backslash Z \to Y \to X] 
$$
for $Z\subset Y$ a closed complex algebraic subvariety of $Y$.
Taking $Z=Y_{red}$ we see that these classes depend only on the underlying
reduced spaces.
By resolution of singularities, $K_{0}(var/X)$ is generated by
classes $[Y\to X]$ with $Y$ smooth, pure dimensional and proper
over $X$. Moreover, for any morphism $f: X'\to X$ we get a functorial pushdown:
$$f_{!}: K_{0}(var/X')\to K_{0}(var/X);\: [h:Z\to X']\mapsto [f\circ h: Z\to X]$$
as well as cross products
$$[Z\to X]\times [Z'\to X'] = [Z\times Z' \to X \times X']\:.$$

By using the motivic Chern class transformation \cite{BSY}
$$mC_y:K_0(var/X) \to G_0(X)[y],$$
with $G_0(X)$ the Grothendieck group of coherent sheaves, 
the un-normalized Hirzebruch class transformation 
$$T_{y*}:=td_* \circ mC_y:K_0(var/X) \to H_*(X)\otimes \bQ[y]$$
was introduced in \cite{BSY} as a class version of the $\chi_y$-genus of a complex algebraic variety $X$. The latter is defined in terms of the Hodge filtration of Deligne's mixed Hodge structure on $H_c^*(X,\bC)$ as:
\be 
\sum_{j,p}\,(-1)^j\dim_{\bC}\Gr_F^p H_c^j(X,\bC)\,(-y)^p =:\chi_y(X)\:.
\ee
Here $td_*:G_0(-)\to H_*(-) \otimes \bQ$ is the Todd class transformation of Baum-Fulton-MacPherson \cite{BFM,Fu2}, with $H_*(-)$ denoting either the Chow group $CH_*(-)$ or the even degree Borel-Moore homology group
$H^{BM}_{2*}(-)$.
The normalized Hirzebruch class transformation, denoted here by $\widehat{T}_{y*}$, provides a 
functorial unification of the (rational) Chern class transformation of MacPherson \cite{MP}, Todd class transformation of Baum-Fulton-MacPherson \cite{BFM} and $L$-class transformation of Cappell-Shaneson \cite{CS1}, respectively. This normalization $\widehat{T}_{y*}$ is obtained by pre-composition the transformation $T_{y*}$ with the normalization functor
$$\Psi_{(1+y)}: H_*(X) \otimes \bQ[y] \to H_*(X) \otimes \bQ[y,(1+y)^{-1}]$$ given in degree $k$ by multiplication by $(1+y)^{-k}$. And it follows from \cite{BSY}[Theorem 3.1] that $\widehat{T}_{y*}:=\Psi_{(1+y)} \circ T_{y*}$ takes in fact values in $H_*(X) \otimes \bQ[y]$, so one is allowed to specialize the parameter $y$ also to the value $y=-1$.

The group homomorphisms $mC_y$, $T_{y*}$ and  $\widehat{T}_{y*}$ are functorial for  push-forwards of proper morphisms, and they commute with cross products. If $X$ is a point, these transformations reduce to the  ring homomorphism 
$$\chi_y:K_0(var/\{pt\}) \to \bZ[y];\: [Z]:=[Z\to \{pt\}]\mapsto \chi_y(Z) \:.$$

The {\it motivic Chern and homology Hirzebruch characteristic classes} of a complex algebraic variety $X$,  denoted by $mC_y(X)$, $T_{y*}(X)$ and resp. $\widehat{T}_{y*}(X)$, are defined as
$$mC_y(X):=mC_y([id_X]) \ ,\quad T_{y*}(X):=T_{y*}([id_X]) \quad \text{and}\quad 
\widehat{T}_{y*}(X):=\widehat{T}_{y*}([id_X])\:.$$
By functoriality, one gets for $X$ compact (so that $X\to \{pt\}$ is proper):
$$\chi_y(X)
=\int_{X}  T_{y*}(X)=\int_{X} \widehat{T}_{y*}(X).$$

Then the following {\it normalization} property holds for the homology Hirzebruch classes of a {\it smooth} complex algebraic variety $X$:
\be\label{norm} T_{y*}(X)=T_y^*(T_X) \cap [X]\:, \ \ \widehat{T}_{y*}(X)=\widehat{T}_y^*(T_X) \cap [X]\:,\ee
with $T_y^*(T_X)$ and $\widehat{T}_y^*(T_X)$ the (un-)normalized  versions of the 
cohomology Hirzebruch class of $X$ appearing in the generalized Riemann-Roch theorem, see \cite{H}[Sect.21.3]. 
More precisely, we have un-normalized and resp. normalized power series
\be\label{ps} Q_y(\alpha):=\frac{\alpha(1+ye^{-\alpha})}{1-e^{-\alpha}} \ ,\,\,\,
\widehat{Q}_y(\alpha):=\frac{\alpha(1+ye^{-\alpha(1+y)})}{1-e^{-\alpha(1+y)}}
\in
\bQ[y][[\alpha]]\ee
with initial terms 
$$Q_y(\al)=1+y+ \frac{1-y}{2} \al + \cdots, \q \widehat{Q}_y(\al)=1+\frac{1-y}{2} \al + \cdots$$
so that for $X$ smooth we have 
\be T^*_y(T_X)=\prod_{i=1}^{\dim X}Q_y(\alpha_i)\in
H^{*}(X) \otimes \bQ[y] \ee
(and similarly for $\widehat{T}_y^*(T_X)$), with $\{\alpha_i\}$  the Chern roots of $T_X$. 
These two power series are related by the following relation
$$ \widehat{Q}_y(\alpha)=(1+y)^{-1}\cdot Q_y(\alpha(1+y)),$$
which explains the use of the normalization functor
$\Psi_{(1+y)}$ in the definition of $\widehat{T}_{y*}(X)$.

Note that for the values $y=-1$, $0$, $1$ of the parameter,  the class $\widehat{T}_y^*$ specializes to the total Chern class $c^*$, Todd class $td^*$, and $L$-class $L^*$, respectively. Indeed, the power series $\widehat{Q}_y(\alpha)\in\bQ[y][[\alpha]]$ becomes respectively
$$1+\alpha,\q\alpha/(1 -e^{-\alpha}),\q\alpha/\tanh(\alpha).$$
On the other hand, the un-normalized class $T_y^*$ specializes for $y=-1$, $0$, $1$ to the top-dimensional Chern class $c^{top}$, Todd class $td^*$, and Atiyah-Singer $L$-class $L_{AS}^*$, respectively, since the power series $Q_y(\alpha)\in\bQ[y][[\alpha]]$ becomes respectively
$$\alpha,\q\alpha/(1 -e^{-\alpha}),\q\alpha/\tanh(\alpha/2).$$

Similarly, we are allowed to specialize the parameter $y$ in the homology classes $ T_{y*}(X)$ and $\widehat{T}_{y*}(X)$
to the distinguished values $y=-1,0,1$. For example, for $y=-1$, we have the identification (cf. \cite{BSY})
\be\label{-1} T_{-1*}(X)=c_0(X) \otimes \bQ \quad \text{and} \quad
\widehat{T}_{-1*}(X)=c_*(X) \otimes \bQ\ee
with the degree zero part $c_0(X)$ and resp. the total (rational) Chern
class $c_*(X)$ of MacPherson \cite{MP}. It is only for this parameter value $y=-1$ that there
is an essential difference between $T_{y*}(X)$ and $\widehat{T}_{y*}(X)$.
Also, for a variety $X$ with at most ``Du Bois singularities"
(e.g., with at most rational singularities \cite{K,Sa}, such as  toric varieties), we have by \cite{BSY} that
\be\label{y0}T_{0*}(X)=\widehat{T}_{0*}(X)=td_*(X):=td_*([\mathcal{O}_X]) \:,\ee
for  $td_*$ the
Baum-Fulton-MacPherson transformation \cite{BFM, Fu2}. And it is still only conjectured that if $X$ is a compact  complex algebraic variety which is also a rational homology manifold
(e.g., a complete simplicial toric variety), then
\be\label{i-35}
\Psi_{2}\left( T_{1*}(X)\right) = \widehat{T}_{1*}(X)=L_*(X)\in H_{2*}(X,\bQ)
\ee 
is the Thom-Milnor homology  $L$-class of $X$ (see \cite{BSY}). For a smooth compact algebraic variety, (\ref{i-35}) follows already from the normalization
property (\ref{norm}).
Moreover, this conjecture is proved in \cite{CMSS}[Cor.1.2] for projective varieties $X$ of the form $Y/G$, with $Y$ a projective manifold and $G$ a finite group of algebraic automorphisms of $Y$.
Hence, the identification (\ref{i-35}) holds for complete algebraic manifolds as well as for projective simplicial toric varieties, which have a finite quotient description by a recent result of
\cite{GS}. 

\vspace{1.5mm}

By definition and functoriality, the motivic Chern and resp. homology Hirzebruch classes satisfy the following additivity property:
\bp\label{p1}
Let $\cS=\{S\}$ be a finite partition of the complex algebraic variety $X$ into disjoint locally closed
complex algebraic subvarieties $S\subset X$, with $i_{\bar{S}}: \bar{S}\to X$ the closed
inclusion of the closure $ \bar{S}$ of $S$. Then
$$[id_X]= \sum_{S}\: [S\hookrightarrow X] = 
\sum_{S}\:(i_{\bar{S}})_{!} [S\hookrightarrow \bar{S}] \in K_0(var/X),$$
so that
\be\label{0} mC_y(X)=\sum_{S}\:mC_y\left([S\hookrightarrow X] \right)
=\sum_{S}\: (i_{\bar{S}})_{*}mC_y\left([S\hookrightarrow  \bar{S}] \right).
\ee
A similar formula holds for the Hirzebruch classes $T_{y*}$ and \tyh.
\ep

Finally, the following crucial result is proved in (greater generality in) \cite{MSS2}[Prop.5.5]
by using an alternative and more sophisticated definition of the Hirzebruch class transformations
in terms of M. Saito's theory of mixed Hodge modules (as explained in \cite{BSY, Sch}):
\bp Let $X$ be a smooth complex algebraic variety, with $D\subset X$ a simple normal crossing divisor, and $i: U:=X\backslash D\hookrightarrow X$ the inclusion of the open complement.
Then 
\be\label{00-1} \begin{split} mC_y \left([i: U\hookrightarrow X]\right)&=
\sum_{p\geq 0} \:[\cO_{X}(-D) \otimes\Omega_X^p(\log D)]\cdot y^p \\
&=[\cO_{X}(-D)
\otimes \Lambda_y \Omega_X^1(\log D)] \in G_0(X)[y],
\end{split}
\ee
and
\be\label{00} \begin{split} T_{y*}\left([i: U\hookrightarrow X]\right)&=
\sum_{p\geq 0} \:td_{*}\left([\cO_{X}(-D) \otimes\Omega_X^p(\log D)]\right)\cdot y^p \\
&=td_{*}\left([\cO_{X}(-D)
\otimes \Lambda_y \Omega_X^1(\log D)]\right) \in H_*(X)\otimes \bQ[y],
\end{split}
\ee
where for a vector bundle $V$ (or for the corresponding locally free sheaf) we set $\Lambda_yV:=\sum_p\Lambda^pV\cdot y^p$, and  in the last formula the transformation $td_*$ is linearly extended
over $\bQ[y]$.
\ep

Alternatively, a new proof of this proposition (without using mixed Hodge modules) will be discussed in the next section.

\subsubsection{Motivic Chern and Hirzebruch classes via the Du Bois complex}
An equivalent definition of the motivic Chern class uses the filtered {\it Du Bois complex} $(\underline{\Omega}^{\bullet}_X,F)$ of $X$, \cite{DB}. More precisely, the following identification holds (see \cite{BSY} for details):
\be\label{dbm} 
mC_y(X)=\sum_{p \geq 0} \ [ \ \uuline{\Omega}^{p}_X] \cdot y^p:=
\sum_{i,p \geq 0} (-1)^i \ [\cH^i(\uuline{\Omega}^{p}_X) ] \cdot y^p \in G_0(X)[y],\ee
where $$\uuline{\Omega}^{p}_X:=\Gr^p_F(\underline{\Omega}^{\bullet}_X)[p] \in D^b_{\rm coh}(X)$$ is a bounded complex of sheaves with coherent cohomology, which coincides with the sheaf of $p$-forms $\Omega^p_X$ on a smooth variety $X$.

For a compact variety $X$, it follows from the degeneration of the corresponding Hodge to De Rham spectral sequence (see \cite{DB}) that:
\be
gr^p_F H^{p+q}(X;\bC)\simeq H^q(X, \uuline{\Omega}^{p}_X), 
\ee 
where $F$ denotes here the Deligne-Hodge filtration on $ H^{p+q}(X;\bC)$. In particular, we then have the identification:
$$\sum_{j,p}\,(-1)^j\dim_{\bC}\Gr_F^p H_c^j(X,\bC)\,(-y)^p =:\chi_y(X) 
=\sum_{p \geq 0} \chi(X,\uuline{\Omega}^{p}_X) \cdot y^p.$$

The un-normalized homology Hirzebruch class of a singular variety $X$  can be given by:
\be\label{db} 
T_{y*}(X)=\sum_{p \geq 0} td_*\left([ \ \uuline{\Omega}^{p}_X] \right) \cdot y^p,\ee
and similarly for the normalized version. 
By definition, $X$ has at most {\it Du Bois singularities} if the canonical map $$\cO_X \to \Gr^0_F(\underline{\Omega}^{\bullet}_X)=\uuline{\Omega}^{0}_X$$ is a quasi-isomorphism, so that in this case we get
$$td_*(X):=td_*([\cO_X])=td_*([\uuline{\Omega}^{0}_X])=:T_{0*}(X).$$
In particular, for $X$ compact:
$$\chi(X,\cO_X)=\chi_0(X).$$
More generally, for a closed algebraic subset $D \overset{i}{\hookrightarrow} X$ of a complex algebraic variety $X$, one has a Du Bois complex $(\underline{\Omega}^{\bullet}_{X,D},F)$ of the pair $(X,D)$ with 
$$\uuline{\Omega}^{p}_{X,D}:=\Gr^p_F(\underline{\Omega}^{\bullet}_{X,D})[p] \in D^b_{\rm coh}(X)$$ 
 fitting for all $p \geq 0$ into an exact triangle (see \cite{K2}[(3.9.2)):
 $$ \uuline{\Omega}^{p}_{X,D} \longrightarrow \uuline{\Omega}^{p}_{X}  \longrightarrow i_* \uuline{\Omega}^{p}_{D} \overset{[1]}{ \longrightarrow} .$$
 In particular, $$[ \ \uuline{\Omega}^{p}_{X,D} ]=[ \ \uuline{\Omega}^{p}_{X}]  -i_*[ \ \uuline{\Omega}^{p}_{D}] \in G_0(X).$$
 Therefore,
 \begin{equation*}
 \begin{split}
 mC_y ( [X \setminus D \hookrightarrow X] ) &= mC_y ( [id_X] ) - i_* mC_y( [id_D] ) \\
 &= \sum_{p \geq 0}  \left(  [ \ \uuline{\Omega}^{p}_X] - i_* [ \ \uuline{\Omega}^{p}_D]  \right) \cdot y^p \\ 
 &=\sum_{p \geq 0}    [ \ \uuline{\Omega}^{p}_{X,D} ]  \cdot y^p.
 \end{split}
 \end{equation*}
Furthermore, if $D$ is a simple normal crossing divisor in a smooth algebraic variety $X$, then there is a filtered quasi-isomorphism (see \cite{K2}[Ex.3.10]):
$$\underline{\Omega}^{\bullet}_{X,D} \simeq \Omega^{\bullet}_D(\log D) \otimes \cO_X(-D),$$ with the trivial filtration on the twisted logarithmic de Rham complex. This reproves formula (\ref{00-1}) from the previous section.

\vspace{1.5mm}

Let us now specialize to the main objects studied in this paper, so let's assume that $X$ is a torus-invariant closed algebraic subset of a toric variety $X_{\Sig}$ (also called a {\it toric polyhedron} in \cite{I}). Then it follows from \cite{I}[Prop.4.2] that there is a natural quasi-isomorphism 
$$
\uuline{\Omega}^{p}_X \simeq \wti{\Omega}^p_X,
$$
with $\wti{\Omega}^p_X$ the {\it Ishida sheaf} of $p$-forms on such a space $X$. In particular, this explains the description (\ref{zero}) of the motivic Chern class and resp. (\ref{abis}) for the un-normalized Hirzebruch class given in the Introduction. Moreover, $X$ has only Du Bois singularities since $\wti{\Omega}^0_X \simeq \cO_X$, see \cite{I}[Cor.4.3]. Note that if $X=X_{\Sig}$ is the ambient toric variety, then the sheaf of Ishida $p$-forms coincides with the sheaf of Zariski $p$-forms on $X_{\Sig}$, see \cite{I}[p.119]. Recall that by definition, the sheaf  of Zariski $p$-forms on a normal variety $X$ is defined as
$$
\wti{\Omega}^p_X:=(j_{\rm reg})_* {\Omega}^p_{X_{\rm reg}},
$$
where $j_{\rm reg}$ is the inclusion of the regular part $X_{\rm reg}$ in $X$.

\br
The identification $
\uuline{\Omega}^{p}_X \simeq \wti{\Omega}^p_X,
$ with the sheaf of Zariski $p$-forms also holds for algebraic $V$-manifolds (also called {\it orbifolds}), i.e., for spaces locally given as the quotient of a complex manifold $Y$ by the action of a finite group $G$ of  automorphisms of $Y$; see \cite{MSS}[Rem.2.9].
In particular, for a {global orbifold} $X=Y/G$, with $Y$ a projective complex algebraic manifold and $G$ a finite group of algebraic automorphisms of $Y$, we have $\wti{\Omega}^p_X \simeq (\pi_*\Omega^p_Y)^G$ for $\pi:Y \to X=Y/G$ the quotient map. So in this case the un-normalized
Hirzebruch class $T_{y*}(X)$ agrees with a corresponding class $$\tau_y(X):=\sum_{p \geq 0} td_* \big( (\pi_*\Omega^p_Y)^G \big) \cdot y^p$$ introduced by Moonen in \cite{Mo}[p.167].
\er


\subsection{A generalized Hirzebruch-Riemann-Roch theorem}

Let $X$ be a compact complex algebraic variety of dimension $d$, with $D$
a Cartier divisor  on $X$. Set 
$$\uuline{\Omega}^{p}_X(D):=\uuline{\Omega}^{p}_X \otimes \cO_X(D) \in D^b_{\rm coh}(X)$$ 
and define the {\it Hirzebruch polynomial of $D$} by the formula
\be
\chi_y(X,\cO_X(D)):=\sum _{p=0}^d \chi(X,\uuline{\Omega}^{p}_X(D)) \cdot y^p.
\ee
Note that in the case of a toric variety $X:=X_{\Sig}$ or an algebraic orbifold $X$,  
$$\uuline{\Omega}^{p}_X(D) \simeq \wti{\om}_X^p \otimes  \cO_X(D)$$
is the sheaf of Zariski $p$-forms twisted by $D$. A similar assertion holds for a torus-invariant closed algebraic subset of a toric variety, by using the sheaf of Ishida $p$-forms instead.

Then we have the following result (which in the case of a global projective orbifold
$X=Y/G$, with $Y$ smooth, is due  to Moonen \cite{Mo}[Satz 2.8, p.169], by using his characteristic class $\tau_y(X)$):
\bt\label{gHRR}(generalized Hirzebruch-Riemann-Roch)\newline
Let $X$ be a compact complex  algebraic varity, and let $D$ be a  Cartier divisor on $X$. Then the Hirzebruch polynomial of $D$ is computed by the formula:
\be\label{Hp}
\chi_y(X,\cO_X(D))=\int_{X} ch(\cO_X(D)) \cap T_{y*}(X).
\ee
\et

\begin{proof}
This is an easy application of the module property of the Todd class transformation: if $\alpha \in K_0(X)$ and $\beta \in K^0(X)$, then the following holds (see \cite{Fu2}[Thm.18.3]):
\be\label{m}
td_*(\beta \otimes \alpha)=ch(\beta) \cap td_*(\alpha),
\ee
with $ch$ denoting the {\it Chern character}.
We therefore have the following sequence of equalities:
\be
\begin{split}
\int_{X} ch(\cO_X(D)) \cap T_{y*}(X)&\overset{(\ref{db})}{=}\sum_{p=0}^d \big[ \int_{X} ch(\cO_X(D)) \cap td_*([ \ \uuline{\Omega}^{p}_X]) \big] \cdot y^p\\
&\overset{(\ref{m})}{=}\sum_{p=0}^d \big[ \int_{X}  td_*([ \ \uuline{\Omega}^{p}_X(D)]) \big] \cdot y^p\\
&\overset{(\ast)}{=}\sum _{p=0}^d \chi(X,\uuline{\Omega}^{p}_X(D)) \cdot y^p\\
&=\chi_y(X,\cO_X(D)),
\end{split}
\ee
where ($\ast$) follows from the degree property of the Todd class transformation.
\end{proof}


\subsection{General facts about toric varieties}\label{basic} In this section we review some basic facts and terminology  from the theory of toric varieties, as needed in this paper. Excellent sources for toric varieties are the books \cite{CLS,Dan,Fu,O}.

Let $M \simeq \bZ^d$ be a $d$-dimensional lattice in $\bR^d$, and let $N=\Hom(M,\bZ)$ be the dual lattice.
Denote the natural pairing by $$\langle \cdot , \cdot \rangle:M \times N \to \bZ.$$
A {\it rational polyhedral cone} $\sig \subset N_{\bR}:=N \otimes \bR$ is a subset of the form 
$$\sig=\con(S):=\{ \sum_{u \in S} \lambda_u \cdot u \ | \  \lambda_u \geq 0 \},$$
for some finite set $S \subset N$. We say that $\sig$ is {\it strongly convex} if $\sig \cap (-\sig)=\{0\}$. The dimension of a cone $\sig$ is the dimension of the subspace of $N_{\bR}$ spanned by $\sig$. Every cone has a {\it dual cone} $\check{\sig}$ defined as
$$\check{\sig}=\{ m \in M_{\bR} \ | \ \langle m,u \rangle \geq 0,  \forall u \in \sig \}.$$
A {\it face} of a cone $\sig$ is a subset $$\tau :=\{ u \in \sig \ | \  \langle m,u \rangle=0 \} \subset \sig $$
for some $m \in M \cap \check{\sig}$. If $\tau$ is a face of $\sig$ we write $\tau \preceq \sig$. Every face of $\sig$ is again a rational polyhedral cone, and a face of a face is a face.

A {\it fan} $\Sig$ in $N_{\bR}$ consists of a finite collection of strongly convex rational polyhedral cones in $N_{\bR}$ so that: 
\bn 
\item if $\sig \in \Sig$, then any face of $\sig$ is also in $\Sig$.
\item if $\sig_1, \sig_2 \in \Sig$, then $\sig_1 \cap \sig_2$ is a face of each (so $\sig_1 \cap \sig_2 \in \Sig$).
\en
The {\it support} of $\Sig$ is the set $|\Sig |:=\bigcup_{\sig \in \Sig} \sig \subset N_{\bR}$. The fan is called {\it complete} if $|\Sig |=N_{\bR}$.   Denote by $\Sigma(i)$ the set of $i$-dimensional cones of $\Sigma$, and similarly, by $\sig(i)$ the collection of $i$-dimensional faces of a cone $\sig$. One-dimensional cones in $\Sig$ are called {\it rays}. For each ray $\rho \in \Sig(1)$, we denote by $u_{\rho}$ the unique generator of the semigroup $\rho \cap N$. The $\{ u_{\rho} \}_{\rho \in \sig(1)}$ are the {\it generators} of $\sig$.  
A fan $\Sig$ is called {\it smooth} if every cone $\sigma \in \Sig$ is smooth, i.e., generated by a part of the $\bZ$-basis of $N$. A fan $\Sig$ is called {\it simplicial} if every cone in $\sigma \in \Sig$ is simplicial, i.e., its generators are linearly independent over $\bR$. The {\it multiplicity} $\mult(\sig)$ of a simplicial cone $\sig \in \Sig$ with generators $u_{\rho_1}, \dots, u_{\rho_k}$ is $| N_{\sig}/{(u_{\rho_1}, \dots, u_{\rho_k})} |$, where $N_{\sig}$ is the sublattice of $N$ spanned by points in $\sig \cap N$, and $|-|$ denotes the cardinality of a (finite) set.

To any fan $\Sig$ one associates a {\it toric variety} $X_{\Sigma}$ as follows. Each cone $\sigma \in \Sig$ gives the affine toric variety $$U_{\sig}=\Sp\big(\bC[M \cap \check \sig]\big)$$
where $\bC[M \cap \check \sig]$ is the $\bC$-algebra with generators $\{ \chi^m \ | \ m \in M \cap \check \sig \}$ and relations $\chi^m \cdot \chi^{m'}=\chi^{m+m'}$. The toric variety $X_{\Sigma}$ is obtained by glueing these affine pieces together, so that $U_{\sig_1}$ and $U_{\sig_2}$ get glued along $U_{\sig_1 \cap \sig_2}$.
The affine toric variety of the trivial cone $\{0\}$ is the torus $T_N=N \otimes \bC^*=\Sp(\bC[M])$, so $T_N$ is an affine open subset of $X_{\Sigma}$. It follows that $M$ can be identified with the character lattice of $T_N$, i.e., $M \simeq \Hom_{\bZ}(T_N,\bC^*)$, while $N\simeq \Hom_{\bZ}(\bC^*,T_N)$ is the lattice of one-parameter subgroups of $T_N$.

The action of the torus $T_N$ on itself extends to an algebraic action of $T_N$ on $X_{\Sigma}$.
Moreover, the $T_N$-action on $X_{\Sigma}$ has only finitely many orbits, and there is an inclusion-reversing correspondence between cones $\sig \in \Sig$ and orbit closures $V_{\sig}:=\overline{O}_{\sig}$. In particular, each ray $\rho \in \Sig(1)$ corresponds to an irreducible $T_N$-invariant divisor $V_{\rho}$ on $X_{\Sigma}$.  A torus-invariant closed algebraic subset $X:=X_{\Sig'}$ of $X_{\Sig}$ (also called {\it toric polyhedron} in \cite{I})  is just a closed union of torus orbits $O_{\sig}$, $\sig \in \Sigma' \subseteq \Sigma$,  with $\Sigma'$ a star-closed subset of $\Sig$ (i.e., if $\tau \in \Sigma'$ is a face of $\sig$, then $\sig \in \Sig'$). 

If $\Sig$ is a fan in $N_{\bR} \simeq \bR^d$, then the toric variety $X_{\Sigma}$ is a complex algebraic variety of dimension $d$ with only mild singularities,
which are rational and Cohen-Macaulay. Moreover, the geometry of $X_{\Sigma}$ is deeply influenced by the properties of the fan $\Sig$. For example:
\bn
\item $X_{\Sigma}$ is {complete} iff $\Sig$ is a complete fan. 
\item $X_{\Sigma}$ contains no torus factor iff $N_{\bR}$ is spanned by the ray generators $u_{\rho}$. 
\item $X_{\Sigma}$ is smooth iff $\Sig$ is a smooth fan. 
\item  $X_{\Sigma}$ is an orbifold iff $\Sig$ is a simplicial fan. 
\en

The toric variety $X_{\Sigma}$ associated to a simplicial fan is called a {\it simplicial toric variety}. Such a variety is a $V$-manifold (orbifold) and therefore also a rational homology manifold, so it satisfies Poincar\'e duality over $\bQ$.  The singular locus of a (simplicial) toric variety $X_{\Sigma}$ is $\bigcup_{\sig \in \Sig_{\rm sing}} V_{\sigma}$, the union being taken over the collection $\Sig_{\rm sing}$ of all singular cones in the fan $\Sigma$. Simplicial toric varieties (without a torus factor) admit a very nice description as geometric quotients \cite{Cox}. This will be explained in some detail in Section \ref{LRR} below.


\subsection{Polytopes and toric varieties}\label{polytop} In this section we explain the connection between polytopes and toric varieties.

A {\it polytope} $P \subset M_{\bR}$ is the convex hull of a finite set $S \subset M_{\bR}$. $P$ is called a {\it lattice polytope} if its vertices lie in $M$. The dimension $\dim({P})$  of a polytope $P$ is the dimension of the smallest affine subspace of $M_{\bR}$ containing $P$. A polytope $P$ whose dimension equals $d=\dim(M_{\bR})$ will be called full-dimensional. Faces of codimension one of $P$ are called {\it facets}.  A polytope $P$ is called {\it simple} if every vertex of $P$ is the intersection of precisely $\dim({P})$ facets. 

From now on, let $P$ be a full-dimensional lattice polytope. Then for each facet $F$ of $P$ there is a unique pair $(u_F,a_F) \in N \times \bZ$ so that $P$ is uniquely described by its {\it facet presentation}:
\be
P=\{ m \in M_{\bR} \ | \ \langle m,u_F \rangle \geq -a_F, \ {\rm for \ all \ facets} \ F \ {\rm of } \ P \}.
\ee
Here $u_F$ is the unique ray generator of the inward-pointing facet normal of $F$.

To a full-dimensional lattice polytope $P$ we associate a fan $\Sig_P$, called the {\it inner normal fan} of $P$, which is defined as follows: to each face $Q$ of $P$ associate the cone $\sig_Q$ by 
$$
\sig_Q:=\con(u_F | F \ {\rm contains } \ Q).
$$
In particular, $\rho_F:=\sig_F=\con(u_F)$ is the ray generated by $u_F$. Let 
$$ X_P:=X_{\Sig_P} $$
 be the associated toric variety. This is commonly referred to as the toric variety associated to the polytope $P$.
As $\Sig_P$ is a complete fan, it follows that $X_P$ is proper.  Moreover, a {\it polytopal subcomplex} $P' \subseteq P$ (i.e., a closed union of faces of $P$) corresponds to a torus-invariant closed algebraic subset $X:=X_{P'}$   of $X_P$. For example, the boundary $\partial P$ of the polytope $P$ corresponds to the union of all toric invariant divisors  $D_F:=D_{\rho_F}$ corresponding to the rays $\rho_F=\sig_F$.

Note that for a cone $\sig_Q \in \Sig_P$ associated to a face $Q$ of $P$, the corresponding orbit closure $V_{\sig_Q}$ can be identified with the toric variety $X_Q$, where $Q$ is the lattice polytope defined as follows: translate $P$ to a vertex of $Q$ so that the origin is a vertex of $Q$; while this operation does not change $\Sigma_P$ or $X_P$, $Q$ is now a full-dimensional polytope in $Span(Q)$. So $Q$ is a full-dimensional lattice polytope relative to $Span(Q) \cap M$, and  $X_Q$ is the associated toric variety.

With the aboth notations, the {\it divisor $D_P$ of the polytope $P$} is defined as:
\be\label{dp}
D_P:=\sum_{F} a_F D_F,
\ee
where $F$ runs over the collection of facets of $P$. Then $D_P$ is a Cartier divisor on $X_P$, and
\be\label{1}
\Gamma(X_P,\cO_{X_P}(D_P)) = \bigoplus_{m \in P \cap M} \bC \cdot \chi^m,
\ee
where $\chi^m$ denotes the character defined by $m \in M$. Moreover, $D_P$ is ample and basepoint free (e.g., see \cite{CLS}[Prop.6.1.10]). Therefore, $X_P$ is a {\it projective} variety.


\subsection{Counting lattice points in lattice polytopes}\label{countsec} Danilov \cite{Dan} used the Hirzebruch-Riemann-Roch theorem to establish a direct connection between the problem of counting the number of lattice points in a lattice polytope and the Todd classes of the associated toric variety (see also \cite{Br, BV2}).

Let $P \subset M_{\bR} \cong \bR^d$ be a full-dimensional lattice polytope with associated projective toric variety $X_P$ and ample Cartier divisor $D_P$ (as described in Section \ref{polytop}). It follows from (\ref{1}) together with the Demazure vanishing theorem \cite{CLS}[Thm.9.2.3] that 
\be\label{2}
\chi(X_P,\cO_{X_P}(D_P))=\dim H^0(X_P,\cO_{X_P}(D_P))= | P \cap M |,
\ee
where $| P \cap M |$ denotes the cardinality of the set $P \cap M$, i.e., the number of lattice points in $P$. 
So the Euler characteristic $\chi(X_P,\cO_{X_P}(D_P))$ computes the number of lattice points in the polytope $P$. Similarly, the Batyrev-Borisov vanishing theorem \cite{CLS}[Thm.9.2.7] yields:
\be\label{3}
\chi(X_P,\cO_{X_P}(-D_P))=(-1)^d \dim H^d(X_P,\cO_{X_P}(-D_P))= 
(-1)^d \cdot | \Int({P}) \cap M |,
\ee
with $\Int({P})$ denoting the interior of $P$.
The left-hand side of equations (\ref{2}) and resp. (\ref{3}) can be  computed by using the Hirzebruch-Riemann-Roch theorem (or the degree property of the Todd class transformation) of \cite{BFM}. For example, we have that
\be\label{4} 
\begin{split}
| P \cap M | =
\chi(X_P,\cO_{X_P}(D_P)) & =\int_{X_P} ch(\cO_{X_P}(D_P)) \cap td_*(X_P),\\
&=\sum_{k\geq 0} \frac{1}{k!} \int_{X_P} [D_P]^k \cap td_k(X_P),
\end{split}
\ee
and similarly
\be\label{4b} 
\begin{split}
| \Int({P}) \cap M | =
(-1)^d\chi(X_P,\cO_{X_P}(-D_P))&=\sum_{k\geq 0} \frac{(-1)^d}{k!} \int_{X_P} [-D_P]^k \cap td_k(X_P)\\
& =\int_{X_P} ch(\cO_{X_P}(D_P)) \cap td_*\left([\omega_{X_P}]\right),
\end{split}
\ee
with $\omega_{X_P}$ the {\it dualizing sheaf} of $X_P$.
Here we use the following duality property of the Todd class transformation $td_*$: if $X$ is a Cohen-Macaulay variety of dim $d$ (e.g., a toric variety $X_{\Sig}$, see \cite{CLS}[Thm.9.2.9]), then for any locally free sheaf $\cE$ on $X$, we have (e.g., see \cite{Fu2}[Ex.18.3.19]):
\be\label{5}
td_k([\cE \otimes \cO_X])=(-1)^{d-k} td_k ([\cE^{\vee} \otimes \omega_X]),
\ee
where $\cE^{\vee}$ denotes the dual of $\cE$.
In particular,
$$td_k(X)=(-1)^{d-k} td_k ([\omega_X]).$$

It thus follows from (\ref{4}) and (\ref{4b}) that counting lattice points in (the interior of) a  full-dimensional lattice polytope $P$ amounts to computing the (dual) Todd class $td_*(X_P)$ (resp. $td_*\left([\omega_{X_P}]\right)$) of the associated projective toric variety. 

Moreover, the  Todd class $td_*(X_P)$ also computes the coefficients of the  {\it Ehrhart polynomial of $P$}. More precisely, if for a positive integer $\ell$, we consider the dilated polytope $\ell P:=\{\ell \cdot u | \ u \in P \}$, the toric varieties corresponding to $P$ and $\ell P$, respectively, are the same, but $D_{\ell P}=\ell D_P$. Therefore, the coefficients of the {\it Ehrhart polynomial} $${\rm Ehr}_P(\ell):=|\ell P \cap M|=\sum_{k=0}^d a_k \ell^k$$ are computed by: $$a_k= \frac{1}{k!}  \int_{X_P}  [D_P]^k \cap td_k(X_P).$$
The duality property (\ref{5}) for the Todd class implies via (\ref{4b}) the {\it Ehrhart reciprocity}
$$(-1)^d \cdot {\rm Ehr}_P(-\ell) =| \Int({\ell P}) \cap M |.$$

In the context of a simplicial toric variety $X_{\Sig}$ (e.g., $X_P$ for $P$ a {\it simple} full-dimensional lattice polytope), a formula for the Todd class $td_*(X_{\Sig})$ was obtained in \cite{BV} via equivariant cohomology, in \cite{EG} by using the Lefschetz-Riemann-Roch theorem, and in \cite{P} by using resolution of singularities and the birational invariance of Todd classes (compare also with \cite{BP, BZ, GP, Mor, P2, P3, PT}).


\section{Characteristic classes of  toric varieties via torus orbit decomposition}\label{dec}

We begin our investigation of motivic Chern and homology Hirzebruch classes of toric varieties by first making use of the torus orbit decomposition together with the additivity properties of these classes as described in Section \ref{motHir}.  This method directly applied to torus-invariant closed algebraic subsets of a toric variety. We also discuss special cases of our computation, yielding formulae for Chern, Todd and resp. $L$-classes of toric varieties.

\subsection{Motivic Chern and Homology Hirzebruch classes of toric varieties}
Let $X_{\Sig}$ be a toric variety of dimension $d$ corresponding to the fan $\Sigma$, and with torus $T=T^d:=(\bC^*)^{d}$. Then $X_{\Sig}$ is stratified by the orbits of the torus action. More precisely, by the orbit-cone correspondence, to a $k$-dimensional cone $\sigma \in \Sig$ there corresponds a $(d-k)$-dimensional torus-orbit \os $ \ \cong (\bC^*)^{d-k}$.  The closure \vs \ of \os \ (which is the same in both classical and Zariski topology) is itself a toric variety, and a $T$-invariant subvariety of $X$. Moreover, if $\sig$ is a face of $\tau$, for short $\sig \preceq \tau$, then $O_{\tau} \subseteq \overline{O}_{\sig}=:V_{\sig}$, and we have that:
$$
V_{\sig}=\bigsqcup_{\sig \preceq \tau} O_{\tau}.
$$
In particular, we can apply the results of Section \ref{motHir} in the setting of toric varieties.  The following preliminary result is a direct consequence of formula (\ref{0}):
\bp\label{p2} Let $X:=X_{\Sig'} \subseteq X_{\Sig}$ be a  torus-invariant closed algebraic subset of $X_{\Sig}$ defined by a star-closed subset $\Sig' \subseteq \Sig$. 
The motivic Chern class of such a space $X$ can be computed by the formula:
\be\label{111} mC_y(X)=\sum_{\sig \in \Sig'}\,mC_y\big([O_{\sig}\hookrightarrow X]\big)=\sum_{\sig \in \Sig'}\,(k_{\sig})_*mC_y\big([O_{\sig}\hookrightarrow V_{\sig}]\big), \ee
with $k_{\sig}:  V_{\sig}   \hookrightarrow  X$ the orbit closure inclusion. A similar formula holds for the homology Hirzebruch classes $T_{y*}(X)$ and  $\widehat{T}_{y*}(X)$ of $X$.
\ep

Moreover, the classes $mC_y\big([O_{\sig}\hookrightarrow V_{\sig}]\big)$ in the above formula can be computed as follows:
\bp\label{p3} For any cone $\sig \in \Sig$, with corresponding orbit \os \ and inclusion $i_{\sig} : O_{\sig} \hookrightarrow \overline{O}_{\sig}=V_{\sig}$, we have:
\be\label{11110}
mC_y\big([O_{\sig}\hookrightarrow V_{\sig}]\big)=(1+y)^{\dim(O_{\sig})} \cdot [\omega_{V_{\sig}}].
\ee
Moreover, 
\be\label{1111}
T_{y*}\big([O_{\sig}\hookrightarrow V_{\sig}]\big)=(1+y)^{\dim(O_{\sig})} \cdot td_*([\omega_{V_{\sig}}]) = (-1-y)^{\dim(O_{\sig})} \cdot 
td_*(V_{\sig})^{\vee}.
\ee
Here $\omega_{V_{\sig}}$ is the canonical sheaf of $V_{\sig}$ and the duality $(-)^{\vee}$
acts by multiplication by $(-1)^i$ on the $i$-th degree homology part $H_i(-)$.\ep

\begin{proof}
The first equality of (\ref{1111}) is a consequence of formula (\ref{00}) and of properties of toric varieties.
More precisely, let $\wti{V}_{\sig} \overset{f_{\sig}}{\to} V_{\sig}$ be a toric resolution of singularities of $V_{\sig}$. Then $\wti{V}_{\sig}$ is a toric variety obtained by refining the fan of $V_{\sig}$. Let $\tilde{i}_{\sig}:O_{\sig} \hookrightarrow \wti{V}_{\sig}$ be the natural open inclusion, with complement the simple normal crossing divisor $D_{\sig}$ whose irreducible components correspond to the rays in the fan of $\wti{V}_{\sig}$. Then $i_{\sig}=f_{\sig} \circ \tilde{i}_{\sig}$, and formula (\ref{00}) applied to the open inclusion  $\tilde{i}_{\sig}:O_{\sig} \hookrightarrow \wti{V}_{\sig}$ yields by the functoriality of $mC_y$ with respect to the proper morphism $f_{\sig}$ that:
\be\label{6} 
mC_y\big([O_{\sig}\hookrightarrow V_{\sig}]\big)=(f_{\sig})_*\big([\cO_{\wti{V}_{\sig}}(-D_{\sig})
\otimes \Lambda_y \Omega_{\wti{V}_{\sig}}^1(\log D_{\sig})]\big).\ee
However, since $\wti{V}_{\sig}$ is a smooth toric variety, the locally free sheaf  $\Omega_{\wti{V}_{\sig}}^1(\log D_{\sig})$ is in fact a trivial sheaf of rank equal to $\dim(\wti{V}_{\sig})=\dim(O_{\sig})$, e.g., see \cite{CLS}[(8.1.5)]. 
So we get that:
\be
[\cO_{\wti{V}_{\sig}}(-D_{\sig})
\otimes \Lambda_y \Omega_{\wti{V}_{\sig}}^1(\log D_{\sig})]=(1+y)^{\dim(O_{\sig})} \cdot  [\cO_{\wti{V}_{\sig}}(-D_{\sig})]  \in G_0(\wti{V}_{\sig})[y].
 \ee

Furthermore, note that the canonical dualizing sheaf $\omega_{\wti{V}_{\sig}}$ on the toric variety $\wti{V}_{\sig}$ is precisely given by $\cO_{\wti{V}_{\sig}}(-D_{\sig})$, e.g, see \cite{CLS}[Thm.8.2.3]. And since the toric morphism $f_{\sig}:\wti{V}_{\sig} \to V_{\sig}$ is induced by a refinement of the fan of $V_{\sig}$, it follows from \cite{CLS}[Thm.8.2.15] that there is an isomorphism:
\be\label{canonical}
{f_{\sig}}_*\omega_{\wti{V}_{\sig}} \simeq \omega_{V_{\sig}}.
\ee
Altogether, we obtain the following sequence of equalities:
\be
\begin{split}
 mC_y\big([O_{\sig}\hookrightarrow V_{\sig}]\big) &=(1+y)^{\dim(O_{\sig})} \cdot (f_{\sig})_*\big( [\cO_{\wti{V}_{\sig}}(-D_{\sig})] \big) \\
 &=(1+y)^{\dim(O_{\sig})} \cdot (f_{\sig})_* \big( [\omega_{\wti{V}_{\sig}}] \big) \\
   &=(1+y)^{\dim(O_{\sig})} \cdot  [\omega_{V_{\sig}}] ,
 \end{split}
\ee
where the last equality is a consequence of (\ref{canonical}) and the vanishing of the higher derived image sheaves of $\omega_{\wti{V}_{\sig}}$, i.e., 
$R^i(f_{\sig})_*\omega_{\wti{V}_{\sig}}=0$, for all $i >0$, 
see \cite{CLS}[Thm.9.3.12].

By applying the Todd class transformation $td_*$ to equation (\ref{11110}), we obtain the first equality of formula (\ref{1111}).
The second equality of (\ref{1111}) follows from the duality property 
(\ref{5}) of the Todd class transformation $td_*$ on the Cohen-Macaulay variety $X_{\Sig}$,
\be\label{duality-td}
td_*([\omega_{V_{\sig}}])= (-1)^{\dim(O_{\sig})} \cdot td_*(V_{\sig})^{\vee}.
\ee
\end{proof}

We can now state the main theorem of this section, which is a direct consequence of Proposition \ref{p2} and Proposition \ref{p3} above:
\bt\label{t1}
Let $X_{\Sig}$ be the toric variety defined by the fan $\Sigma$, and $X:=X_{\Sig'} \subseteq X_{\Sig}$ a  torus-invariant closed algebraic subset defined by a star-closed subset $\Sig' \subseteq \Sig$. For each cone $\sig \in \Sig$ with corresponding orbit $O_{\sig}$, denote by $k_{\sig} : V_{\sig} \hookrightarrow X$ the inclusion of the orbit closure.   Then the 
motivic Chern class $mC_y(X)$ is computed by the formula:
\be\label{i7m} 
mC_y(X)=\sum_{\sig \in \Sig'}\, (1+y)^{\dim(O_{\sig})} \cdot (k_{\sig})_*[\omega_{V_{\sig}}].
\ee
Similarly, the un-normalized homology Hirzebruch class ${T}_{y*}(X)$ is computed by:
\be\label{7} 
\begin{split}
T_{y*}(X)&=\sum_{\sig \in \Sig'}\, (1+y)^{\dim(O_{\sig})} \cdot (k_{\sig})_*td_*([\omega_{V_{\sig}}])\\
&=\sum_{\sig \in \Sig'}\, (-1-y)^{\dim(O_{\sig})} \cdot 
(k_{\sig})_*\left( td_*(V_{\sig})^{\vee}\right).
\end{split}
\ee
\et

If $X$ in the above theorem is compact (i.e., $\Sig$ is a complete fan), then by taking the degree in formula (\ref{7}) we get as a corollary the well-known formula for the Hodge polynomial of toric varieties and resp. compact toric polyhedra:
\bc Let $X:=X_{\Sig'} \subseteq X_{\Sig}$ be a torus-invariant compact algebraic subset of the toric variety $X_{\Sig}$. Then:
$$
\chi_y(X) =\sum _{\sig \in \Sig'}\, (-1-y)^{\dim(O_{\sig})}.
$$
\ec 

We conclude this section by giving a formula for the normalized Hirzebruch class $\widehat{T}_{y*}(X)$ of such a toric polyhedron $X$. Recall that \tyh \ is defined by applying the renormalization $\Psi_{(1+y)}$ to $T_{y*}$, i.e., by multiplying each component in $H_k(-)\otimes \bQ[y]$ by $(1+y)^{-k}$. We therefore obtain the following result:

\bt\label{t2}
Let $X_{\Sig}$ be the toric variety defined by the fan $\Sigma$, and $X:=X_{\Sig'} \subseteq X_{\Sig}$ a  torus-invariant closed algebraic subset defined by a star-closed subset $\Sig' \subseteq \Sig$. For each cone $\sig \in \Sig$ with corresponding orbit $O_{\sig}$, denote by $k_{\sig} : V_{\sig} \hookrightarrow X$ the inclusion of the orbit closure.   Then the normalized homology Hirzebruch class $\widehat{T}_{y*}(X)$ is computed by the following formula:
\be\label{8} 
\begin{split}
\widehat{T}_{y*}(X)&=\sum_{\sig\in \Sig', k}\, (1+y)^{\dim(O_{\sig}) -k} \cdot (k_{\sig})_*td_k([\omega_{V_{\sig}}])\\
&=\sum_{\sig \in \Sig', k}\, (-1-y)^{\dim(O_{\sig})-k}  \cdot (k_{\sig})_*td_k(V_{\sig}).
\end{split}
\ee
\et


\subsection{Chern, Todd, and $L$-classes of toric varieties}\label{ctl}
In this section, we discuss special cases of our formulae (\ref{7}) and (\ref{8}) for $X:=X_{\Sig'} \subseteq X_{\Sig}$ a  torus-invariant closed algebraic subset defined by a star-closed subset $\Sig' \subseteq \Sig$, which are obtained by specializing the parameter $y$ to the distinguished values $y=-1, 0, 1$.
For each cone $\sig \in \Sig$ with corresponding orbit $O_{\sig}$, denote as before by $k_{\sig} : V_{\sig} \hookrightarrow X$ the inclusion of the orbit closure. 

First, recall that by making $y=-1$ in the normalized Hirzebruch class $\widehat{T}_{y*}$ one gets the (rational) MacPherson Chern class $c_*$. So we get as a corollary of Theorem \ref{t2} the following result from \cite{BBF}:
\bc[Ehler's formula]
The (rational) MacPherson-Chern class $c_*(X)$ is computed by:
\be\label{E}
c_*(X)=\sum_{\sig \in \Sig'}\,  (k_{\sig})_*td_{\dim(O_{\sig})}(V_{\sig})=\sum_{\sig \in \Sig'}\,  (k_{\sig})_*[V_{\sig}].
\ee
\ec

Ehler's formula for the MacPherson-Chern class $c_*(X_{\Sig})$ (even with integer coefficients) is due to
\cite{Al,BBF}. In fact, Theorems \ref{t1} and \ref{t2} should be regarded as counterparts of Ehler's formula for the motivic Chern and Hirzebruch classes. Moreover, our proof of Theorem \ref{t1} follows closely the strategy of Aluffi \cite{Al} for his proof of  Ehler's formula.\\

Secondly, since $X$ is DuBois, we get for $y=0$ the following formula for the Todd class:
\bc  The Todd class $td_*(X)$ of $X$ is computed by:
\be\label{T}
td_*(X)=\sum_{\sig \in \Sig'}\, (k_{\sig})_*td_*([\omega_{V_{\sig}}])
=\sum_{\sig \in \Sig'}\, (-1)^{\dim(O_{\sig})} \cdot 
(k_{\sig})_*\left( td_*(V_{\sig})^{\vee}\right).
\ee
\ec

By applying the duality (\ref{duality-td}) to (\ref{T}) in the case $X=X_{\Sig}$ (which is Cohen-Macauley), one also gets in the above notations the following
\be\label{3dual}
td_*([\omega_{X_{\Sig}}])=
\sum_{\sig \in \Sig}\, (-1)^{{\rm codim}(O_{\sig})} \cdot 
(k_{\sig})_* td_*(V_{\sig})\:.
\ee

Finally, we have the following result about the classes  $T_{1*}(X)$
and  $\widehat{T}_{1*}(X)$ of a toric polyhedron $X$, respectively the homology $L$-classes of compact toric manifolds and projective simplicial toric varieties $X_{\Sig}$:
\bc The classes $T_{1*}(X)$
and  $\widehat{T}_{1*}(X)$ are computed by:
\be\label{4un} 
T_{1*}(X)= \sum_{\sig \in \Sig'}\, 2^{\dim(O_{\sig})} \cdot (k_{\sig})_*td_*([\omega_{V_{\sig}}]) =
\sum_{\sig\in \Sig'}\, (-2)^{\dim(O_{\sig})} \cdot 
(k_{\sig})_*\left( td_*(V_{\sig})^{\vee}\right)
\ee
and
\be\label{4nor} 
\widehat{T}_{1*}(X)=\sum_{\sig\in \Sig', i}\, (-2)^{\dim(O_{\sig})-i}  \cdot (k_{\sig})_*td_i(V_{\sig}).
\ee
Moreover, $T_{1*}(X_{\Sig})=L^{AS}_*(X_{\Sig})$ is the Atiyah-Singer $L$-class for 
$X_{\Sig}$ a compact toric manifold, and $\widehat{T}_{1*}(X_{\Sig})=L(X_{\Sig})$
is the Thom-Milnor $L$-class for $X_{\Sig}$ a compact toric manifold
or a simplicial projective toric variety.
\ec

\begin{proof} This follows by making $y=1$ in formulae (\ref{7}) and (\ref{8}), after noting that the following identifications hold:
$$
L^{AS}_*(X_{\Sig})=T_{y*}(X_{\Sig})|_{y=1} \quad \text{and} \quad
L_*(X_{\Sig})=\widehat{T}_{y*}(X_{\Sig})|_{y=1}.
$$
For the latter identity, we use \cite{CMSS}[Cor.1.2] together with the fact that a simplicial projective toric variety $X_{\Sig}$ can be realized globally as a {\it finite} quotient (see \cite{GS}). 
\end{proof}

By taking degrees in formula (\ref{4nor}), we get as a consequence the following signature formula
(compare also with \cite{LR}):

\bc The signature of a compact toric manifold or a projective simplicial toric variety $X_{\Sig}$ can be computed as:
$$
sign(X_{\Sig})=\chi_1(X_{\Sig})= \sum_{\sig \in \Sig}\, (-2)^{\dim(O_{\sig})}.
$$
\ec

Note that the  equality $sign(X)=\chi_1(X)$ is also a consequence of the {\it Hodge index theorem} for a projective orbifold $X$, which is a special case of the corresponding result
in intersection homology for a projective variety $X$ (see \cite{MSS}[Sec.3.6], and in the smooth case see also \cite{O}[Thm.3.12]).\\

Let us finish this section with the following duality result:
\bp\label{4Dual}  Let $X_{\Sig}$ be a simplicial toric variety of dimension $d$. Then
\be\label{4dual} 
T_{1*}(X_{\Sig})= (-1)^{d}\cdot \left(T_{1*}(X_{\Sig})\right)^{\vee}\:.
\ee
\ep

\begin{proof}
Note that this result is trivial for a toric manifold by the normalization
property (\ref{norm}), since the corresponding characteristic class $T_{1}^*=L^*_{AS}$
of the tangent bundle is defined by the even power series $\alpha/\tanh(\alpha/2)$.
Since the homology Hirzebruch class $T_{1*}$ also commutes with cross products (\cite{BSY}),
we can also assume that the simplicial toric variety $X_{\Sig}$ contains no torus factor.
Then we would like to have a duality result similar to (\ref{duality-td})  for all sheaves of Zariski
$p$-forms $\wti{\Omega}^p_X$ coming in the description of $T_{1*}$  in (\ref{abis}):
$$T_{1*}(X_{\Sig})=\sum_{p\geq 0}\; td_*\left([\wti{\Omega}^p_X]\right).$$
Note that, while the coherent sheaves $\wti{\Omega}^p_X$ are in general not locally free, 
in the case of a simplicial toric variety $X_{\Sig}$ they are maximal Cohen-Macauly by \cite{CLS}[p.418]. Hence,
$$\cE xt^q_{\cO_{X_{\Sig}}}(\wti{\Omega}^p_X, \omega_{X_{\Sig}})=0 \quad 
\text{for $q>0$.}$$
Finally, we have (e.g., see \cite{CLS}[Ex.8.0.13])
$$\cH om_{\cO_{X_{\Sig}}}(\wti{\Omega}^p_X, \omega_{X_{\Sig}})\simeq 
\wti{\Omega}^{d-p}_X ,$$
so that
$$td_*\left([\wti{\Omega}^p_X]\right) =( -1)^{d}\cdot td_*\left([\wti{\Omega}^{d-p}_X]\right)^{\vee}$$
for all $0\leq p \leq d$, as in \cite{Fu2}[Ex.18.3.19].
\end{proof}

\br Note that the duality property (\ref{4dual}) holds more generally for any complex algebraic variety $X$ of pure dimension $d$, which, moreover, is a rational homology manifold, so by Verdier duality we get that: $D(\bQ_X)=\bQ_X[2d](d)$. Then the result follows from this equality on the mixed Hodge module level by \cite{Sch}[Rem.5.20]. \er

Applying (\ref{4dual}) to (\ref{4un}), one also gets for a simplicial toric variety $X_{\Sig}$ the following equality:
\be\label{simpl.dual} 
T_{1*}(X_{\Sig})= (-1)^{\dim(X_{\Sig})}\cdot \left(T_{1*}(X_{\Sig})\right)^{\vee}= 
\sum_{\sig\in \Sig}\, (-1)^{\dim(X_{\Sig})}\cdot (-2)^{\dim(O_{\sig})} \cdot 
(k_{\sig})_* td_*(V_{\sig}).
\ee

\br Our formulae in this section calculate the homology Hirzebruch classes (and, in particular, the Chern, Todd and $L$-classes, respectively) of singular toric varieties in terms of the Todd classes of closures of torus orbits. Therefore, known results about Todd classes of toric varietes (e.g., see  \cite{BP, BeV, BV, BZ, GP, Mor, P, P2, P3, PT}) also apply to these other characteristic classes: in particular, Todd class formulae satisfying the so-called {\it Danilov condition}, as well as the algorithmic computations from \cite{BP, P3}.
In the simplicial context, these  Todd classes can also be computed by using the Lefschetz-Riemann-Roch theorem of \cite{EG} for geometric quotients. This will be recalled in Section \ref{LRR}.
\er


\section{Weighted lattice point counting via Hirzebruch classes}
In this section we show that, just like the Todd classes of toric varieties can be used for counting the number of lattice points in a lattice polytope (see Section \ref{countsec}), the un-normalized homology Hirzebruch classes yield  lattice point counting with certain weights. \\

Let $P \subset M_{\bR} \cong \bR^d$ be a full-dimensional lattice polytope with associated projective toric variety $X_P$ and ample Cartier divisor $D_P$ (as described in Section \ref{polytop}). Recall from (\ref{4}) and (\ref{4b}) that  the number of lattice points in (the interior of) $P$ can be  computed by the formulae:
$$ 
| P \cap M | =
\chi(X_P,\cO_{X_P}(D_P))=\int_{X_P} ch(\cO_{X_P}(D_P)) \cap td_*(X_P),
$$
and 
$$ 
| \Int({P}) \cap M | =
(-1)^d\chi(X_P,\cO_{X_P}(-D_P))=\int_{X_P} ch(\cO_{X_P}(D_P)) \cap td_*\left([\omega_{X_P}]\right).
$$

In what follows, we use the notation $Q \preceq P$ for a face $Q$ of the polytope $P$. Note that
$$
| P \cap M |=\sum_{Q \preceq P} | \Relint(Q) \cap M |,
$$
where the number of points inside a face $Q$  of $P$ is counted with respect to the lattice $Span(Q) \cap M$. Combining these formulae, we therefore get:
\be\label{comb}
\chi(X_P,\cO_{X_P}(D_P))=\sum_{Q \preceq P} | \Relint(Q) \cap M |.
\ee

As a generalization of  (\ref{comb}), the main result of this section is a weighted lattice point counting formula for lattice polytopes and for their closed subcomplexes,
 which motivates our calculation of homology Hirzebruch classes via the torus-orbit decomposition (cf. Theorem \ref{t1}):
 
\bt\label{wc} Let $M$ be a lattice of rank $d$ and $P \subset M_{\bR} \cong \bR^d$ be a full-dimensional lattice polytope with associated projective toric variety $X_P$ and ample Cartier divisor $D_P$. In addition, let $X:=X_{P'}$ be a torus-invariant closed algebraic subset of $X_P$ corresponding to a polytopal subcomplex $P' \subseteq P$ (i.e., a closed union of faces of $P$). Then the following formula holds:
\be\label{wcf}
\sum_{Q \preceq P'} (1+y)^{\dim(Q)} \cdot | \Relint(Q) \cap M |=\int_{X} ch(\cO_{X_P}(D_P)|_X) \cap T_{y*}(X),
\ee
where the summation on the left is over the faces $Q$ of $P'$, $|-|$ denotes the cardinality of sets and $\Relint(Q)$ the relative interior of a face $Q$.
\et

\begin{proof} For a face $Q$ of $P'$, denote by $i_Q:V(\sig_Q):=X_Q \hookrightarrow X$ the inclusion of the orbit closure associated to the (cone of the) face $Q$. Note that we have $\dim(Q)=\dim(O_{\sig_Q}).$ Then by (\ref{7}) the following  equality holds:
$$\int_{X_P} ch(\cO_{X_P}(D_P)|_X) \cap T_{y*}(X)
= \sum_{Q \preceq P'} (1+y)^{\dim(Q)} \int_{X} ch(\cO_{X_P}(D_P)|_X) \cap (i_Q)_*td_*([\omega_{X_Q}]).$$

It remains to prove that for any face $Q$ of $P'$, we have that:
\be\label{lem2}
 \int_{X} ch(\cO_{X_P}(D_P)|_X) \cap (i_Q)_*td_*([\omega_{X_Q}]) = | \Relint(Q) \cap M |.
\ee
This follows from the functorial properties of the cap product. Indeed,
\begin{equation*}
\begin{split}
 \int_{X} ch(\cO_{X_P}(D_P)|_X) \cap (i_Q)_*td_*([\omega_{X_Q}]) 
 &= \int_{X_Q} (i_Q)^* ch(\cO_{X_P}(D_P)|_X) \cap td_*([\omega_{X_Q}]) \\
 &= \int_{X_Q}  ch((i_Q)^*(\cO_{X_P}(D_P)|_X)) \cap td_*([\omega_{X_Q}]) \\
&= \int_{X_Q}  ch(\cO_{X_Q}(D_Q)) \cap td_*([\omega_{X_Q}]) \\
&\overset{(\ref{4b})}{=}| \Relint(Q) \cap M |.
 \end{split}
\end{equation*}
\end{proof}

Recall that, if $\ell$ a positive integer, and we consider the dilated polytope $\ell P:=\{\ell \cdot u | \ u \in P \}$ and its subcomplex $\ell P'$, then the toric varieties (resp. toric polyhedra) corresponding to $P$, $\ell P$ and $P'$, $\ell P'$,  respectively, are the same, but $D_{\ell P}=\ell D_P$. Therefore, we have the following:

\bc
The Ehrhart polynomial of the polytopal subcomplex $P' \subseteq P$ as above, counting the number of lattice points in the dilated complex $\ell P':=\{\ell \cdot u \ | \ u \in P' \}$ for a positive integer $\ell$, is computed by:
$${\rm Ehr}_{P'}(\ell):=|\ell P' \cap M|=\sum_{k=0}^{d'} a_k \ell^k,$$ with $$a_k= \frac{1}{k!}  \int_{X}  [D_P|_X]^k \cap td_k(X),$$ and $X:=X_{P'}$ the associated toric polyhedron of dimension $d'$.
\ec

Note that ${\rm Ehr}_{P'}(0)=a_0=\int_X td_*(X)$ can be computed on the combinatorial side by formula (\ref{y}) as the Euler characteristic of $P'$, while on the algebraic geometric side it is given by the arithmetic genus of $X:=X_{P'}$: 
\be\label{Er} \chi({P'}):=\sum_{Q \preceq P'} (-1)^{\dim Q}=\chi_0(X)=\int_X td_*(X)=\chi(X,\cO_X).\ee

\vspace{1.5mm}

As another application of Theorem \ref{wc}, by specializing at $y=1$, 
we get for a full-dimensional lattice polytope $P$ the following identity:
\be\label{wcL}
\sum_{Q \preceq P} \left(\frac{1}{2}\right)^{{\rm codim}(Q)} \cdot |\Relint(Q) \cap M |=\int_{X_P} ch(\cO_{X_P}(D)) \cap \left( \frac{1}{2}\right)^{\dim(X_P)}\cdot T_{1*}(X_P) .
\ee
 Moreover, by
(\ref{simpl.dual}) one gets in the same way also the following result for a full-dimensional {\it simple} lattice polytope $P$:
\be\label{5L}
\sum_{Q \preceq P} \left(-\frac{1}{2}\right)^{{\rm codim}(Q)} \cdot |Q \cap M |=\int_{X_P} ch(\cO_{X_P}(D)) \cap \left( \frac{1}{2}\right)^{\dim(X_P)}\cdot T_{1*}(X_P).
\ee

Note that if $X_P$ is smooth (i.e., $P$ is a {\it Delzant} lattice polytope), the un-normalized
Hirzebruch class $T_{1*}(X_P)$ corresponds under Poincar\'e duality to the characteristic
class of the tangent bundle defined in terms of Chern roots by the un-normalized power series
$\alpha/\tanh(\frac{\alpha}{2})$, so that $\left( \frac{1}{2}\right)^{\dim(X_P)}\cdot T_{1*}(X_P)$
corresponds to the characteristic
class of the tangent bundle defined in terms of Chern roots  by the normalized power series
$\frac{\alpha}{2}/\tanh(\frac{\alpha}{2})$. 

Therefore (\ref{5L}) extends the result 
\cite{Fe}[Corollary 5] for Delzant lattice polytopes to the more general case of simple lattice polytopes. Our approach is very different from the techniques used in \cite{Fe} for projective toric manifolds, which rely on virtual submanifolds Poincar\'e dual to cobordism characteristic classes of the projective toric manifold.\\

Together with the generalized Hirzebruch-Riemann-Roch
Theorem \ref{gHRR} one also gets from Theorem \ref{wc} the following combinatorial description:

\bc\label{chiy-M} Let $M$ be a lattice of rank $d$ and $P \subset M_{\bR} \cong \bR^d$ be a full-dimensional lattice polytope with associated projective simplicial toric variety $X=X_P$ and ample Cartier divisor $D_P$. In addition, let $X:=X_{P'}$ be a torus-invariant closed algebraic subset of $X_P$ corresponding to a polytopal subcomplex $P' \subseteq P$ (i.e., a closed union of faces of $P$). Then the following formula holds:
\be
\chi_y(X,\cO_X(D_P|_X)) = \sum_{Q \preceq P'} (1+y)^{\dim(Q)} \cdot | \Relint(Q) \cap M |.
\ee
\ec

\br Let $d':=\dim(P')$ and denote by $P'(i)$ the set of $i$-dimensional faces of $P'$. Then the result of Corollary \ref{chiy-M} is equivalent to the following formula:
\be\label{Mat}
\chi(X,\wti{\om}^p_{X} (D_P|_X))=\sum_{i=p}^{d'}  {i \choose p} \sum_{Q \in P'(i)} | \Relint(Q) \cap M |,
\ee
with $\wti{\om}^p_{X} $ the sheaf of Ishida $p$-forms on $X$.
Indeed, 

\be
\begin{split}
\chi_y(X,\cO_{X}(D_P|_X))&:=\sum _{p=0}^{d'} \chi(X,\wti{\om}^p_{X}(D_P|_X)) \cdot y^p\\
&\overset{(\ref{Mat})}{=}\sum _{p=0}^{d'} \big( \sum_{i=p}^{d'}  {i \choose p} \sum_{Q \in P'(i)} | \Int(Q) \cap M | \big) \cdot y^p\\
&=\sum_{i=0}^{d'} \big( \sum_{p=0}^i {i \choose p} y^p  \big) \sum_{Q \in P'(i)} | \Relint(Q) \cap M | \\
&=\sum_{i=0}^{d'} (1+y)^i \sum_{Q \in P'(i)} | \Relint(Q) \cap M | \\
&=\sum_{Q \preceq P'} (1+y)^{\dim(Q)} \cdot | \Relint(Q) \cap M |.
\end{split}
\ee
In the case when $P'=P$ is a simple full-dimensional polytope, formula (\ref{Mat}) was already proved by Materov \cite{M} (see also \cite{CLS}[Thm.9.4.7, Ex.9.4.14]).
\er


\section{Characteristic classes of simplicial toric varieties via Lefschetz-Riemann-Roch}\label{LRR}
In this section we first recall the Cox construction of simplicial toric varieties as geometric quotients, see \cite{Cox} and \cite{CLS}[Sec.5.1]. We then use the Lefschetz-Riemann-Roch theorem of \cite{EG} to compute the homology Hirzebruch classes  of such varieties. We will use freely the notations from Section \ref{basic}.

\subsection{Simplicial toric varieties as geometric quotients}\label{Cox}

Let $\Sig$ be a fan in $N_{\bR} \simeq \bR^d$, with associated toric variety $X_{\Sig}$ and torus $T:=T_N =(\bC^*)^d$. For each ray $\rho \in \Sig(1)$, denote by $u_{\rho}$ the unique generator of the semigroup $\rho \cap N$. Recall that $N$ is identified with the one-parameter subgroups of $T$, i.e., $N \simeq \Hom_{\bZ}(\bC^*,T)$. Let $n=| \Sig(1) |$ be the number of rays in the fan $\Sig$. First we assume that $X_{\Sig}$ contains {\it no torus factor}, i.e. that the ray generators $u_{\rho}$ ($\rho \in \Sig(1)$) span $N_{\bR}$ (so that $n>0$).
This is for example always the case if $\Sig$ is complete.

Define the map of tori $$\gamma:(\bC^*)^n \lra T \quad \text{by} \quad (t_{\rho})_{\rho} \mapsto \prod_{\rho \in \Sig(1)} u_{\rho}(t_{\rho}),$$ 
and let $G:=\ker (\gamma)$. Then $G$ is a product of a torus and a finite abelian group. In particular, $G$ is reductive. Let $Z(\Sig) \subset \bC^n$ be the subvariety defined by the monoidal ideal $$B(\Sig):=\langle  \hat{x}_{\sig} : \sig \in \Sig \rangle$$
where $$\hat{x}_{\sig}:=\prod_{\rho \notin \sig(1)} x_{\rho}$$ 
and $(x_{\rho})_{\rho \in \Sig(1)}$ are coordinates on $\bC^n$. Then the variety 
$$W:=\bC^n \setminus Z(\Sig)$$
is a toric manifold, and there is a toric morphism $$\pi:W \to X_{\Sig}.$$ Indeed, let $\{e_{\rho}  | \ \rho \in \Sig(1) \}$  denote the standard basis of the lattice $\bZ^{n}$. For each cone $\sig \in \Sig$, define the cone \be\label{up} \wti{\sig}:=\con\big(e_{\rho}  | \ \rho \in \sig(1)\big) \subseteq \bR^{n}.\ee
Then these cones and their faces form a smooth fan $\wti{\Sig}$ in $\bR^{n}$ whose associated toric variety is exactly $W$. Moreover, the map of lattices $\bZ^{n} \to N$ sending $e_{\rho} \mapsto u_{\rho}$ induces the toric morphism $\pi:W \to X_{\Sig}$, and for any $\sigma \in \Sig$ the restriction of $\pi$ over the affine piece $U_{\sig}$ is precisely the map $\pi_{\sig}:=\pi|_{U_{\wti{\sig}}}: U_{\wti{\sig}} \to U_{\sig}$. 

The group $G$ acts on $W$ by the restriction of the diagonal action of $(\bC^*)^{n}$, and the toric morphism $\pi$ is constant on $G$-orbits. Then Cox \cite{Cox} proved the following
(see also \cite{CLS}[Thm.5.1.11]): \bt
Let $X_{\Sig}$ be a toric variety without torus factors. Then $X_{\Sig}$ is simplicial if and only if $X_{\Sig}$ is a geometric quotient $W/G$.\et

Let us assume now that $X_{\Sig}$ is such a simplicial toric variety without torus factors. Under the quotient map $\pi:W \to X_{\Sig}$, the coordinate hyperplane $\{ x_{\rho}=0 \}$ maps to the invariant Weil divisor $V_{\rho}$. More generally, if $\sig$ is a $k$-dimensional cone generated by rays $\rho_1,\dots,\rho_k$, then the orbit closure $V_{\sig}$ is the image under $\pi$ of the linear subspace $W_{\sig}$ defined by the ideal $(x_{\rho_1},\dots,x_{\rho_k})$ of the total coordinate ring $\bC\big[ x_{\rho} \ | \ \rho \in \Sig(1) \big]$. Let $G_{\sig}$ be the stabilizer of $W_{\sig}$. Then, 
if we denote by $$a_{\rho}:(\bC^*)^{n} \to \bC^*$$ the projection onto the $\rho$-th factor (and similarly for any restriction of this projection),  we have by definition that
\be\label{stab}
\begin{split}
G_{\sig}&=\{g \in G \ | \ a_{\rho}(g)=1, \forall \rho \notin \sig(1) \}\\
&\simeq \{(t_{\rho})_{\rho \in \sig(1)}  \ | \ t_{\rho}\in \bC^*, \prod_{\rho \in \sig(1)} u_{\rho}(t_{\rho})=1 \}.
\end{split}
\ee

Note that the second description of $G_{\sig}$ in (\ref{stab}) depends only on $\sig$
(and not on the fan $\Sig$ nor the group $G$).
For a $k$-dimensional rational simplicial cone $\sig$ generated by the rays $\rho_1,\dots,\rho_k$ one gets by \cite{CLS}[Prop.1.3.18]
the finite abelian group 
\be\label{Gsig-int}
G_{\sig}\simeq N_{\sig}/(u_1,\dots,u_k)
\ee
defined  as the quotient of the sublattice $N_{\sig}$ of $N$ spanned by the points in $\sig\cap N$ modulo
 the sublattice $(u_1,\dots,u_k)$ generated by the ray generators $u_j$ of  $\rho_j$.
 Or, in other words,  $G_{\sig}=ker(\gamma_{\sig})$, with 
 $$\gamma_{\sig}:  (\bC^*)^k\to T_{N_{\sig}}=(\bC^*)^k,$$ 
 is the group appearing in the quotient construction of
 the toric variety $X_{\sig}\simeq \bC^k/G_{\sig}$ associated the simplicial cone $\sig \subset (N_{\sig})_{\bR}$. Geometrically, $X_{\sig}$ is a ``normal slice'' to the orbit
 $$(\bC^*)^{d-k}\simeq O_{\sig}\subset U_{\sig}\simeq X_{\sig}\times (\bC^*)^{d-k}.$$
Then $|G_{\sig}|=\mult(\sig)$ is just the multiplicity of $\sig$, with $\mult(\sig)=1$ exactly in the case of a smooth cone.
 Let $m_i\in M_{\sig}$ for $1\leq i\leq k$ be the unique primitive elements in the dual lattice $M_{\sig}$ of $N_{\sig}$ satisfying $ \langle m_i,u_j \rangle = 0$ for $i\neq j$ and
 $ \langle m_i,u_i \rangle > 0$
 so that the dual lattice $M_{\sig}'$ of $(u_1,\dots,u_k)$ is generated by the elements $\frac{m_j}{\langle m_j,u_j \rangle}.$
  Then, for $g=n+(u_1,\dots,u_k)\in G_{\sig}=N_{\sig}/(u_1,\dots,u_k)$, the character $a_{\rho_j}(g)$ is also given by (see \cite{Fu}[p.34]):
 \be\label{a-intr}
 a_{\rho_j}(g)=\exp\left(2\pi i\cdot \gamma_{\rho_j}(g)
 \right), \quad \text{with} \quad \gamma_{\rho_j}(g):= 
 \frac{\langle m_j,n \rangle}{\langle m_j,u_j \rangle}\:.  
 \ee


\subsection{Lefschetz-Riemann-Roch theorem for simplicial toric varieties}
In this section, we recall from \cite{EG} the statement of the Lefschetz-Riemann-Roch theorem for simplicial toric varieties.

Let $X=X_{\Sig}$ be a $d$-dimensional simplicial toric variety without a torus factor. In the notations of the previous section, $X$ is the geometric quotient $W/G$ coming from the Cox construction. Denote by $K^G(W)$ the Grothendieck group of $G$-equivariant locally free sheaves on $W$. This is a ring, with product given by the tensor product of equivariant vector bundles. Let $K_G(W)$ be the Grothendieck group of $G$-equivariant coherent sheaves; this is a $K^G(W)$-module. Since $W$ is smooth and $G$ is reductive, we have that $K^G(W)\simeq K_G(W)$. If we identify the representation ring $R(G)$ with $K^G(pt)$, then the pullback by the constant map $W \to pt$ makes $K^G(W)$ into an $R(G)$-algebra and $K_G(W)$ into an $R(G)$-module.
Moreover, if $R:=R(G) \otimes \bC$ denotes the complex representation ring of $G$, then $R$ coincides with the coordinate ring of $G$.

Let $CH^*_G(W) \otimes \bQ$ denote the equivariant rational Chow ring of $W$ in the sense of Edidin-Graham, and denote by $\phi$ the following isomorphism:
$$\phi:CH^*_G(W) \otimes \bQ \overset{\simeq}{\lra} CH^*(X) \otimes \bQ= CH_{d-*}(X) \otimes \bQ$$
or its composition with the cycle map 
$$CH_{d-*}(X) \otimes \bQ \lra H^{\BM}_{2d-2*}(X;\bQ).$$
For $\rho \in \Sig(1)$, let $z_{\rho}:=[W_{\rho}]_G \in CH^1_G(W) \otimes \bQ$ be the equivariant fundamental class of the invariant 
divisor $\{x_{\rho}=0\}$ on $W$. Then we have that 
$$\phi(z_{\rho})=[V_{\rho}] \in H^1(X)\otimes \bQ \simeq H_{d-1}(X) \otimes \bQ,$$ with $H^*(X)$ denoting either the  Chow ring $CH^*(X)$ or the even degree  cohomology ring $H^{2*}(X;\bZ)$.
In the following we usually think of $[V_{\rho}]$ as a cohomology class.
More generally, for any cone $\sig \in \Sig$, it follows that $\phi([W_{\sig}]_G)=\frac{1}{\mult(\sig)} [V_{\sig}]$.

In the notations of the previous section, consider the finite set $$G_{\Sig}:=\bigcup_{\sig \in \Sig} G_{\sig} \subset G=\Sp({R}),$$ and note that the smooth cones of $\Sig$ only contribute the identity element.
For $g \in G$, the set of coordinates of $\bC^{n}$ (with $n=|\Sig(1)|$) on which $g$ acts non-trivially is identified with: $$g(1)=\{ \rho \in \Sig(1) \ | \ a_{\rho}(g) \neq 1 \}.$$ Thus the fixed locus of $g$ in $\bC^n$ is the linear subspace defined by the ideal 
$\big(x_{\rho} : \rho \in g(1)\big)$ of the total coordinate ring. So the fixed locus $W^g$ of $g$ in $W$ is the intersection of this linear subspace with $W=\bC^n \setminus Z(\Sig)$. Note that $W^g \neq \emptyset$ if and only if the linear space defined by the ideal $\big(x_{\rho} : \rho \in g(1)\big)$ is not contained in $Z(\Sig)$, and the latter is equivalent to the existence of a cone $\sig \in \Sig$ which is generated by the rays of $g(1)$. Note that for such a cone $\sig$, we have that $g \in G_{\sig}$ and $W^g=W_{\sig}$.

\vspace{2mm}

Assume  that $X_{\Sig}$ is a $d$-dimensional simplicial toric variety containing a torus factor, and let $N'$ be the intersection of $N$ with the proper subspace of $N_{\bR}$ spanned by all ray generators $u_{\rho}$. Then all cones of $\Sig$ lie in $N'_{\bR}$, so they form a simplicial fan $\Sig'$ in $N'_{\bR}$, with
$$X_{\Sig}\simeq X_{\Sig',N'}\times (\bC^*)^r$$
so that $X_{\Sig',N'}= W'/G'$ contains no torus factor (assuming that $\Sig$ is not just the zero cone $\{0\}$). 
Letting $G:=G'$ act {\it trivially} on $(\bC^*)^r$, we get (non-canonically):
$$X_{\Sig}\simeq X_{\Sig',N'}\times (\bC^*)^r = W/G,$$
with $W:=W'\times (\bC^*)^r$. 
If $\Sig$ is just the zero cone $\{0\}$, we have $X_{\Sig}=(\bC^*)^d=W/G$ with
$W:=(\bC^*)^d$ and $G:=\{id\}$.
So the preceeding discussion and notations can also be extended to this context.

\vspace{2mm}

Let $\ch^G$ and $\Td^G$ denote the equivariant Chern character and resp. equivariant Todd class in equivariant rational Chow cohomology (as in \cite{EG}). We then have the following: 
\bt\label{t-lrr} (Lefschetz-Riemann-Roch theorem for simplicial toric varieties \cite{EG}) \newline
Let $X=X_{\Sig}$ be a simplicial toric variety, and let $\pi:W \to X=W/G$ be the quotient map as before. If $g \in G$, let $i_g: W^g \hookrightarrow W$ denote the inclusion of the fixed locus of $g$, with normal bundle $N_{W^g}$,  
and let  $T_{W^g}$ be the tangent bundle of $W^g$. If $\cF \in K^G(W)$, then:
\be\label{lrr}
td_*\big( [(\pi_* \cF)^G] \big)=\sum_{g \in G_{\Sig}} \phi \circ( i_g)_* \left( \frac{\ch^G\big(i_g^*\cF\big)(g)}{\ch^G(\La_{-1}\big(N^*_{W^g})(g)\big)} \Td^G(T_{W^g}) 
\right),
\ee
with $(i_g)_*$ the equivariant Gysin map induced by Poincar\'e duality.
\et
Here, the (equivariant) {\it twisted Chern character} $\ch^G(-)(g)$ for $g\in G$ is defined 
for a (virtual) $G$-equivariant vector bundle (or locally free sheaf) $\cF \in K^G(W^g)$ as follows:
$$\ch^G(\cF)(g):= \sum_{\chi}\: \chi(g)\cdot \ch^G(\cF_{\chi}),$$
where $\cF\simeq \oplus_{\chi}\: \cF_{\chi}$ is the (finite) decomposition of $\cF$
into subbundles (or sheaves) $\cF_{\chi}$ on which $g$ acts by a character $\chi: \langle g \rangle \to \bC^*$. The  sheaves $\cF_{\chi}$ are also $G$-equivariant, since $G$ is 
{\it abelian}. Note that since these characters are complex valued, the individual terms on the right-hand side of (\ref{lrr})  are in $H_*(X) \otimes \bC$. However, their sum is in $H_*(X) \otimes \bQ$, since this is the case for the left-hand side of (\ref{lrr}). 

\vspace{2mm}

As a corollary of Theorem \ref{t-lrr}, Edidin-Graham proved the following formula for the Todd classes of simplicial toric varieties (obtained for $\cF=[\cO_W]$). In the complete case, the same result follows from the equivariant Riemann-Roch theorem of Brion-Vergne \cite{BV} (compare also with \cite{G} for a related approach using Kawasaki's orbifold Riemann-Roch theorem, with applications to Euler-MacLaurin formulae).
\bc\label{Tc-EG} The Todd class of simplicial toric variety $X=X_{\Sig}$ is computed by the formula:
\be\label{toddc}
td_*(X):=td_*([\cO_X])=\left( \sum_{g \in G_{\Sig}}  \prod_{\rho \in \Sig(1)} \frac{[V_{\rho}]}{1-a_{\rho}(g)e^{-[V_{\rho}]}} \right)
\cap [X]\:.
\ee
\ec


\subsection{Homology Hirzebruch classes of simplicial toric varieties}\label{hirlrr}
The main result of this section is the following computation of homology Hirzebruch classes via the Lefschetz-Riemann-Roch theorem:
\bt\label{Tc} Let $X=X_{\Sig}$ be a simplicial toric variety of dimension $d$, with $n:=| \Sig(1) |$.
The un-normalized Hirzebruch class of  $X=X_{\Sig}$ is computed by the formula:
\be\label{Hirzc}
 T_{y*}(X)= (1+y)^{d-n}\cdot\left(\sum_{g \in G_{\Sig}}  \prod_{\rho \in \Sig(1)} \frac{[V_{\rho}] \cdot 
\big( 1+y  \cdot a_{\rho}(g)  \cdot e^{-[V_{\rho}]}\big)}{1-a_{\rho}(g) \cdot e^{-[V_{\rho}]}} \right) \cap [X]\:,
\ee
and the normalized Hirzebruch class of  $X=X_{\Sig}$ is given by:
\be\label{nHirzc}
\widehat{T}_{y*}(X)= \left( \sum_{g \in G_{\Sig}}  \prod_{\rho \in \Sig(1)} \frac{[V_{\rho}] \cdot 
\big( 1+y  \cdot a_{\rho}(g)  \cdot e^{-[V_{\rho}](1+y)}\big)}{1-a_{\rho}(g) \cdot e^{-[V_{\rho}](1+y)}} \right) \cap [X]\:.
\ee
\et

\begin{proof} The proof follows the lines of  \cite{EG}[Thm.4.2], with some technical additions dictated by our more general setting. If $\Sig$ is just the zero cone $\{0\}$ (so that $n=0$),
with $X=X_{\Sig}=(\bC^*)^d$, then the identities 
$$T_{y*}(X)= (1+y)^{d}\cdot [X] \quad \text{and} \quad \widehat{T}_{y*}(X)=[X]$$
 follow from the normalization condition, since $T_X$ is a trivial vector bundle of rank $d$.
 Moreover, by the multiplicativity of $T_{y*}$ and $\widehat{T}_{y*}$ with respect to cross products, we can assume in the following that $X_{\Sig}$ contains no torus factor. 

Recall from formula (\ref{abis}) that the un-normalized Hirzebruch class of a  (simplicial) toric variety, 
can be computed as:
\be\label{simp}
T_{y*}(X)=\sum_{p=0}^d td_*(\wti{\Omega}^p_X) \cdot y^p=td_*(\big[ \Lambda_y \wti{\Omega}^1_X \big]),
\ee
with $\wti{\Omega}^p_X$ the sheaf of Zariski $p$-forms on $X$ and $ \Lambda_y \wti{\Omega}^1_X:= \sum_{p=0}^d \wti{\Omega}^p_X\cdot y^p$.
We claim that the following identity holds:
\be\label{tech}
\La_y\wti{\Omega}^1_X \simeq \big(\pi_*\La_y N^*_{(G,W)} \big)^G.
\ee
Here $N^*_{(G,W)}$ is by definition the (sheaf of sections of the) dual bundle of the $G$-equivariant quotient bundle 
$$N_{(G,W)}:= T_W/T_{(G,W)},$$  with  $T_{(G,W)}$ the sub-bundle of
$T_W$ given by the vectors tangent to the $G$-orbits. So $N^*_{(G,W)}$ is the sub-bundle (or subsheaf) of $\Omega^1_W$ given by the one-forms 
vanishing along the $G$-orbits. If $a: G\times W\to W$ denotes the action $(g,w)\mapsto g\cdot w$, then the differential
$da: (T_G)_{id}\times (W\times \{0\})\to T_W$ is injective, since $G$ acts with finite stabilizers, and the image of $da$  gets identified with
$T_{(G,W)}$. So $T_{(G,W)}\simeq (T_G)_{id}\times W$, with the adjoint $G$-action on $(T_G)_{id}$, is a trivial $G$-vector bundle of rank
$n-d$ since $G$ is abelian.

So we want to show that for any $p$ we have:
\be\label{tech2}
\wti{\Omega}^p_X \simeq \big(\pi_*\Lambda^p N^*_{(G,W)} \big)^G.
\ee
To prove (\ref{tech2}), let $n=| \Sig(1) |$  and cover $W$ by principal opens $$\{x \in \bC^n : \hat{x}_{\sig} \neq 0 \}$$ associated to the top-dimensional cones $\sig \in \Sig(d)$ of the fan $\Sig$. Then in the notation of (\ref{up}), we have the identification 
$$ \{x \in \bC^n : \hat{x}_{\sig} \neq 0 \}\simeq  \bC^d \times (\bC^*)^{n-d}  \simeq      U_{\wti{\sig}},$$ where the last equality follows  since ${\wti{\sig}}$ is a smooth $d$-dimensional cone in $\bR^n$. Moreover, $G_{\sig}$ acts naturally on $\bC^d$ with:
$$
U_{\wti{\sig}} /G \simeq \bC^d/{G_{\sig}} \simeq U_{\sig}. 
$$
Since proving (\ref{tech2}) is a local question (in fact the following canonical identifications glue together), we can assume
$W=U_{\wti{\sig}} $, $X=U_{\sig}$, $G=G_{\sig} \times  (\bC^*)^{n-d}$, and 
$\pi=\pi_{\sig}$. So $\pi$ can be factored as a free quotient followed by a finite quotient:
$$U_{\wti{\sig}}  \stackrel{\tilde{\pi}}{\lra} U_{\wti{\sig}} /{ (\bC^*)^{n-d}}\simeq \bC^d \lra \bC^d/{G_{\sig}} \simeq U_{\sig}.$$
And in these charts we have the identification $T_{(G,W)}=T_{\tilde{\pi}}$ between $T_{(G,W)}$ and the sub-bundle $T_{\tilde{\pi}}$ of vectors tangent to the fibers of $\tilde{\pi}$. Hence 
$N^*_{(G,W)}=\tilde{\pi}^*\Omega^1_{\bC^d}$.
Then (\ref{tech2}) follows from the known fact 
$$\Omega^p_X \simeq \big(\tilde{\pi}_*\tilde{\pi}^*\Omega^p_X \big)^G$$
about a free  quotient map of algebraic manifolds, and resp. the corresponding decription 
of the sheaf of Zariski $p$-forms 
$$\wti{\Omega}^p_X \simeq \big(\pi_*\Omega^p_M \big)^G$$
for a finite quotient map $\pi: M\to X=M/G$ with $M$ smooth.

Putting (\ref{simp}) and (\ref{tech2}) together, and using the Lefschetz-Riemann-Roch formula (\ref{lrr}), we get:
\be\label{Hlrr}
\begin{split}
T_{y*}(X) &=td_*\big( [ \big(\pi_*\La_y N^*_{(G,W)} \big)^G] \big) \\
& =\sum_{g \in G_{\Sig}} \phi \circ( i_g)_* 
\left( \frac{\ch^G\big(i_g^*\La_y N^*_{(G,W)}(g)\big)}{\ch^G(\La_{-1}\big(N^*_{W^g})(g)\big)} \Td^G(T_{W^g}) 
\right),
\end{split}
\ee
with $i_g:W^g \hookrightarrow W$ the inclusion of the fixed point locus of $g \in G_{\Sig}$, and 
$$\phi:CH^*_G(W) \otimes \bQ \lra H^*(X)\times \bQ \simeq H_{d-*}(X) \otimes \bQ$$ 
is as in (\ref{lrr}). Here, $H^*(X)$ denotes either the Chow ring or the even degreel cohomology ring, and $H_*(X)$ is either the Chow group or the even degree Borel-Moore homology. Recall that for $\rho \in \Sig(1)$, we let $z_{\rho}:=[W_{\rho}]_G \in CH^1_G(W)\otimes \bQ$ denote the equivariant fundamental class of the invariant divisor $\{x_{\rho}=0\}$ on $W$, and we have that $\phi(z_{\rho})=[V_{\rho}]$.

Since $$T_W \simeq \bigoplus_{\rho \in \Sig(1)} \cO_W(x_{\rho}),$$
it follows that the equivariant normal bundle $N_{W^g}$ of $W^g \subset W$ splits as
 $$N_{W^g} \simeq \bigoplus_{\rho \in g(1)} i_g^*\cO_W(x_{\rho}).$$
 This is also a decomposition into eigenbundles for the finite subgroup $\langle g \rangle$, where the weight on the $\rho$-th factor is the character $a_{\rho}$. Since $c_1(\cO_W(x_{\rho}))=z_{\rho}$, we therefore get:
 $$
 \ch^G(\La_{-1}\big(N^*_{W^g})(g)\big)=i_g^* \left(  \prod_{\rho \in g(1)} \big(1-a_{\rho}(g^{-1})e^{-z_{\rho}}\big) \right).
 $$
 Similarly, the tangent bundle of $W^g$ splits as 
 $$T_{W^g} \simeq \bigoplus_{\rho \notin g(1)} i_g^*\cO_W(x_{\rho}).$$
 Hence,
 $$
 \Td^G(T_{W^g}) =i_g^* \left(  \prod_{\rho \notin g(1)} \frac{z_{\rho}}{1-e^{-z_{\rho}}} \right)=i_g^* \left(  \prod_{\rho \notin g(1)} \frac{z_{\rho}}{1-a_{\rho}(g^{-1})e^{-z_{\rho}}} \right),
 $$
 where in the last equality we use the fact that $a_{\rho}(g)=1$ if $\rho \notin g(1)$. From the short exact sequence of $G$-vector bundles
 $$0\to N^*_{(G,W)} \to T^*_W \to T^*_{(G,W)}\to 0$$
 we get 
 $$ \ch^G\big(i_g^*\La_y N^*_{(G,W)}(g)\big) = (1+y)^{d-n}\cdot  \ch^G\big(i_g^*\La_yT^*_W(g)\big)\:,$$
 since $T^*_{(G,W)}$ is a trivial $G$-vector bundle of rank $n-d$.
 Finally,

 $$
 \ch^G\big(i_g^*\La_yT^*_W(g)\big)=\prod_{\rho \in \Sig(1)} i_g^*\ch^G \big( 1+ y \cdot \cO_W(x_{\rho})^{\vee}(g)  \big) =i_g^* \left( \prod_{\rho \in \Sig(1)} \big(  1+ya_{\rho}(g^{-1})e^{-z_{\rho}} \big) \right).
 $$
 
Since $W^g$ is the intersection of $W$ with the linear subspace cut out by the ideal generated by $\{x_{\rho}\}_{\rho \in g(1)}$, we have:
$$ ( i_g)_*( i_g)^*(-)=\big( \prod_{\rho \in g(1)} z_{\rho} \big) \cdot (-).$$
 
 Altogether, we obtain:
 \be
 T_{y*}(X)=(1+y)^{d-n}\cdot \sum_{g \in G_{\Sig}} \phi \left(   \prod_{\rho \in \Sig(1)} \frac{z_{\rho}(1+ya_{\rho}(g^{-1})e^{-z_{\rho}})}{1-a_{\rho}(g^{-1})e^{-z_{\rho}}}  
   \right).
 \ee
 The first result about the un-normalized Hirzebruch class follows,
 since $\phi(z_{\rho})=[V_{\rho}]$ and we use the fact that $g \in G_{\Sig} \iff g^{-1} \in G_{\Sig}$.
 To get the formula for the normalized Hirzebruch class, we just have to substitute $(1+y)^{-d}\cdot [X]$ for $[X]\in H_{d}(X) \otimes \bQ$ and
 $(1+y)\cdot [V_{\rho}]$ for $[V_{\rho}]\in H^1(X)\otimes \bQ$.
\end{proof}

As a by-product of formula (\ref{tech2}) in the proof of Theorem \ref{Tc} above, we also get the following presentation of the motivic Chern class $mC_y(X_{\Sig})$:
\bc In the above notations, the motivic Chern class $mC_y(X_{\Sig})$ of the simplicial toric variety $X_{\Sig}=W/G$ is given by:
\be
mC_y(X_{\Sig})=\pi^G_*\left(  \Lambda_y([T^*_W]\ominus [T^*_{(G,W)}]) \right) \in G_0(X_{\Sig})[y],
\ee
with $\pi^G_*:K^0_G(W) \to G_0(X_{\Sig})$ induced by the corresponding exact functor $\cF \mapsto (\pi_*\cF)^G$.
\ec

As another consequence of Theorem \ref{Tc}, we get the following formula for the Todd class $td_*([\omega_X])$ of the canonical bundle, which played an essential role in Section \ref{dec}.
\bc If $X=X_{\Sig}$ is a $d$-dimensional simplicial toric variety, then
\be\label{toddcan}
td_*([\omega_X])=\left(\sum_{g \in G_{\Sig}}  \prod_{\rho \in \Sig(1)} \frac{ 
 a_{\rho}(g) \cdot [V_{\rho}] \cdot e^{-[V_{\rho}]}}{1-a_{\rho}(g) \cdot e^{-[V_{\rho}]}} \right) \cap [X]\:.
\ee
\ec
\begin{proof}
Recall from formula (\ref{abis}) that the un-normalized Hirzebruch class of a  $d$-dimensional (simplicial) toric variety, can be computed as:
$$
T_{y*}(X)=\sum_{p=0}^d td_*(\wti{\Omega}^p_X) \cdot y^p.
$$ 
So formula (\ref{toddcan}) is obtained by specializing our formula (\ref{Hirzc}) to the top-degree $d$ in $y$.
\end{proof}

Note that if $\Sig$ is a smooth fan, then $G_{\Sig}=\{id_G\}$. In particular, we derive from (\ref{nHirzc}) the following formula for the 
normalized Hirzebruch classes of smooth toric varieties:
\bc\label{sm}
The normalized Hirzebruch class of a smooth toric variety $X=X_{\Sig}$ is given by:
\be\label{Hirzsm}
\widehat{T}_{y*}(X)= \left( \prod_{\rho \in \Sig(1)} \frac{[V_{\rho}] \cdot \big( 1+y \cdot e^{-[V_{\rho}](1+y)}\big)}{1-e^{-[V_{\rho}](1+y)}} \right)
\cap [X] =
\left( \prod_{\rho \in \Sig(1)} \widehat{Q}_y([V_{\rho}] )\right) \cap [X]\:,
\ee
where $\widehat{Q}_y(-)$ is the power series from (\ref{ps}).
\ec

\br By making $y=-1, 0, 1$ in (\ref{nHirzc}) and (\ref{Hirzsm}) we obtain formulae for the Chern, Todd and $L$-classes of simplicial and resp. smooth toric varieties.
\er

\bex\label{fake} {\it (Fake) Weighted projective spaces} \\
Consider a complete simplicial toric variety $X_{\Sig}$ of dimension $d$, whose fan $\Sigma$ has $d+1$ rays $\rho_i$ with primitive generators $u_i$ ($i=1,\dots,d$),  so that $N_{\bR}=\sum_{i=0}^d \bR_{\geq 0} u_i$. Then there exist positive integers $q_0,\dots,q_d$ with $\gcd(q_0,\dots ,q_d)=1$, so that $$q_0 u_0+\dots +q_d u_d=0.$$
By the proof of \cite{RT}[Thm.2.1], these weights are already reduced, i.e., $$1=\gcd(q_0,\dots , q_{j-1}, q_{j+1},\dots q_d),$$ for all $j=0,\dots,d$, so the toric variety $X_{\Sig,N'}$ associated to the fan $\Sig$ and the lattice $N'$ generated by the vectors $u_i$ is the  weighted projective space $\bP(q_0,\dots,q_d)$; see loc.cit. In the case when $N'\neq N$, $X_{\Sig}:=X_{\Sig,N}$ is the {\it fake weighted projective space} defined as the quotient of  $X_{\Sig,N'}=\bP(q_0,\dots,q_d)$ by the finite abelian group $N/N'$, see also \cite{Bu,Ka,RT}.

Then by the Cox construction, we have that $X_{\Sig}=W/G$, with $W=\bC^{d+1}\setminus \{0\}$ and $G\simeq \bC^* \times H$, for a finite abelian group $H$. Note that $\bC^*$ acts on $W$ with weights $(q_0,\dots ,q_d)$, e.g., see \cite{RT}[Thm.1.26], so that $W/{\bC^*}=\bP(q_0,\dots,q_d)$, and $H \simeq N/N'$. By \cite{CLS}[Thm.12.4.1], the rational cohomology ring $H^*(X_{\Sig})\otimes \bQ$ of $X_{\Sig}$ is given as the quotient of the polynomial ring $\bQ[x_0,\dots,x_d]$ by the relations $x_0 \cdots  x_d=0$ and $q_jx_i=q_ix_j$ (for all $i\neq j$), under the identification $x_i=[V_i]$ corresponding to the ray $\rho_i$. Therefore, $H^*(X_{\Sig}) \otimes \bQ \simeq \bQ[T]/(T^{d+1})$ with distinguished generator $T=\frac{[V_{j}]}{q_j}\in H^2(X_{\Sig};\bQ)$, $j=0,\dots, d$. 

The normalized Hirzebruch class of $X_{\Sig}$ is then computed by:
\be\label{wph00}
\widehat{T}_{y*}(X_{\Sig})= \left( \sum_{g \in G_{\Sig}} \: \prod_{j=0}^d \frac{q_j T \cdot 
\big( 1+y  \cdot a_j(g)  \cdot e^{-q_j (1+y) T}\big)}{1-a_j(g)  \cdot e^{-q_j (1+y) T}} \right) \cap [X_{\Sig}]\:.
\ee
In particular, if $N=N'$, i.e., for the weighted projective space $\bP(q_0,\dots,q_d)$, we obtain:
\be\label{wph01}
\widehat{T}_{y*}(X_{\Sig})= \left( \sum_{\lambda \in G_{\Sig}} \: \prod_{j=0}^d \frac{q_j T \cdot 
\big( 1+y  \cdot \lambda^{q_j}  \cdot e^{-q_j (1+y) T}\big)}{1-\lambda^{q_j}  \cdot e^{-q_j (1+y) T}} \right) \cap [X_{\Sig}]\:,
\ee
with $$G_{\Sig}:=\bigcup_{j=0}^d \{ \lambda\in \bC^* \ | \ \lambda^{q_j}=1 \}.$$ 
\eex

\br
The formulae of the above example can be used for obtaining a closed expression for the Ehrhart polynomial of a lattice {\it $d$-simplex} in $\bR^d$, similar to those proved in \cite{C,DR} by very different methods. Indeed, the toric variety associated to such a simplex is a fake weighted projective space. For example, let us consider the simplex $P=conv(0,a_1e_1,\dots, a_de_d)$ with all $a_i \in \bZ_{>0}$, $\gcd(a_1,\dots,a_d)=1$ and $\{e_1,\dots,e_d\}$ the standard basis of $\bZ^d$. By \cite{RT}[Algorithm 3.17], the associated toric variety is the weighted projective space $X_{P}=\bP(1,q_1,\dots,q_d)$ with weights $q_i=\frac{s_i}{s}$, for $s_i:=\frac{\prod_{j=1}^d a_j}{a_i}$ and $s=\gcd(s_1,\dots,s_d)$. Moreover, in this case, the class of the ample Cartier divisor $D_P$ is given by ${\rm lcm}(q_1,\dots,q_d) \cdot T$ (see \cite{RT}[Prop.1.22]), with $T^d \cap [X_P]=\frac{1}{q_1 \cdots q_d}[pt]$.
\er
\subsection{Mock characteristic classes of simplicial toric varieties}\label{M}

We have seen in the previous section that the normalized Hirzebruch class of a smooth toric variety $X=X_{\Sig}$ is computed by the formula
(\ref{Hirzsm}) of Corollary \ref{sm}.
Following the terminology from \cite{CS, P}, we define the (normalized) {\it mock Hirzebruch class} $\widehat{T}^{(m)}_{y*}(X)$ of a 
{\it simplicial} toric variety $X=X_{\Sigma}$ to be the class computed using the formula in the smooth case: 
\be\label{mock} \widehat{T}^{(m)}_{y*}(X):= \left(\prod_{\rho \in \Sig(1)} \frac{[V_{\rho}] \cdot 
\big( 1+y   \cdot e^{-[V_{\rho}](1+y)}\big)}{1-e^{-[V_{\rho}](1+y)}} \right) \cap [X] =
\left( \prod_{\rho \in \Sig(1)} \widehat{Q}_y([V_{\rho}] )\right) \cap [X] \:.\ee
In general, let us pretend that the stable complex bundle of $X_{\Sigma}$ is a sum of line bundles with Chern classes  
$[V_{\rho}] \in H_{n-1}(X) \otimes \bQ \simeq H^1(X) \otimes \bQ$, $\rho \in \Sigma(1)$, and take the appropriate polynomial in these, e.g., 
as in \cite{H}[\S 12] for smooth varieties. In particular, besides mock Hirzebruch  classes, we can define mock Todd classes (see \cite{P}), 
or mock Chern and $L$-classes, respectively. Note that the mock Chern, mock Todd, and mock $L$-classes, respectively, are all specializations 
of the (normalized) mock Hirzebruch class, obtained by evaluating the parameter $y$ at $-1$, $0$, and $1$, respectively. 

\br Let $cl^*$ be a multiplicative characteristic class of vector bundles defined in terms of Chern roots $\alpha$ by a power series $f(\alpha)$. Then the definition of the corresponding
mock characteristic class of a simplicial toric variety $X$, i.e., 
$$cl^{(m)}_*(X):=
\left(\prod_{\rho \in \Sig(1)} \: f([V_\rho])\right)
\cap [X]$$
  agrees in the case of a toric manifold $X$ with the class
$cl^*(T_X)\cap [X]$ only if $f$ is a {\it normalized} power series with  $f(0)=1$
(as can be seen easily by using the exact sequence from \cite{CLS}[Thm.8.1.6]).
So in the case of the Atiyah-Singer $L$-class $L_{AS}^*$ corresponding to the un-normalized
power series $f(\alpha)=\alpha/\tanh(\frac{\alpha}{2})$ with $f(0)=\frac{1}{2}$, it is more useful to define the corresponding mock class in terms of the normalized power series $(\frac{\alpha}{2})/\tanh(\frac{\alpha}{2})$.
For a smooth toric manifold $X$ of dimension $d$, this corresponds then to the identity
$$\left(\frac{1}{2}\right)^d\cdot L^{AS}_*(X)=\left(\frac{1}{2}\right)^d\cdot T_{1*}(X)$$ 
(compare also with the formulae (\ref{wcL}) and
(\ref{5L})).
\er

The main result of this section is a reformulation of our Theorem \ref{Tc} in terms of mock characteristic classes of closures of torus orbits. Before stating the result, let us recall our notations and terminology from Section \ref{Cox}. A simplicial toric variety $X=X_{\Sigma}$
without a torus factor can be realized as a geometric quotient $X=W/G$, with $W \subset \bC^n$ (for $n=|\Sig(1)|$) and $G$ a reductive abelian group. If $\sig$ is a $k$-dimensional cone generated by rays $\rho_1,\dots,\rho_k$, then the orbit closure $V_{\sig}$ is the image under $\pi:W \to X$ of the linear subspace $W_{\sig}$ cut out by the ideal $(x_{\rho_1},\dots,x_{\rho_k})$ of the total coordinate ring $\bC\big[ x_{\rho} \ | \ \rho \in \Sig(1) \big]$. The stabilizer $G_{\sig}$ of $W_{\sig}$
is given by 
$$G_{\sig}=\{g \in G \ | \ a_{\rho}(g)=1, \forall \rho \notin \sig(1) \},$$
with $a_{\rho}:(\bC^*)^{n} \to \bC^*$ the projection onto the $\rho$-th factor. $G_{\sig}$ is a finite abelian group of order equal to the multiplicity of the cone $\sig$. For any $g \in G$, we use as before the notation $$g(1)=\{ \rho \in \Sig(1) \ | \ a_{\rho}(g) \neq 1 \}.$$ Then the fixed locus $W^g$ of $g$ in $W$ is nonempty if and only if $g \in G_{\sig}$ for some cone $\sig \in \Sig$, and if this is the case then $W^g=W_{\sig}$.

In what follows we write as usual $\tau \preceq \sig$ for a face $\tau$ of $\sig$, and we use the notation $\tau \prec \sig$ for a proper face $\tau$ of $\sig$. Note that
$$\tau \preceq \sig  \Longrightarrow G_{\tau} \subseteq G_{\sig}.$$
Set $$G_{\sig}^{\circ}:=G_{\sig} \setminus \bigcup_{\tau \prec \sig} G_{\tau}=\{ g \in G_{\sig} \ | \ a_{\rho}(g) \neq 1, \forall \rho \in \sig(1) \}.$$
Since $G_{\sig}=\{id_G\}$ for a smooth cone $\sig$, it follows that $G_{\sig}^{\circ}=\emptyset$ if $\sig$ is a smooth cone of positive dimension, while $G_{\{0\}}^{\circ}=G_{\{0\}}=\{id_G\}$.
Moreover, as $$G_{\sig}=\bigsqcup_{\tau \preceq \sig}G_{\tau}^{\circ},$$ (with $\bigsqcup$ denoting disjoint unions) 
it follows that \be\label{last} G_{\Sig}:=\bigcup_{\sig \in \Sig} G_{\sig}=\bigsqcup_{\sig \in \Sig} G_{\sig}^{\circ}=\{id_G\} \sqcup \bigsqcup_{\sig \in \Sig_{\rm sing}} G_{\sig}^{\circ},
\ee
where $\Sig_{\rm sing}$ denotes the collection of the singular cones in $\Sig$.
Moreover, in the above notations, the following easy result holds:
\bl\label{l10} For any $\sig \in \Sig_{\rm sing}$ the following are equivalent:
\bn
\item[(a)] $g \in G_{\sig}^{\circ}$.
\item[(b)] $g(1)=\sig(1) $.
\item[({c})] $W^g=W_{\sig}$.
\en
\el
\begin{proof} It follows by definition that ($b$) and (${c}$) are equivalent, so it suffices to show that ($a$) and (${b}$) are equivalent. 

If $g \in G_{\sig}^{\circ}$, then $g \in G_{\sig}$, so $a_{\rho}(g)=1$ for all $\rho \notin \sig(1)$. Hence $g(1) \subseteq \sig(1)$. If moreover, there exists $\rho_1 \in \sig(1) \setminus g(1)$, let $\tau$ denote the face of $\sig$ generated by rays $\rho \in \sig(1) \setminus \{\rho_1 \}$. Since $a_{\rho_1}(g)=1$, it follows that $g \in G_{\tau}$, which is a contradiction since $\tau$ is a proper face of $\sig$. So we must have $g(1)=\sig(1) $.

Conversely, if $g(1)=\sig(1)$, then the cone $\sig$ is generated by the rays in $g(1)$, so $g \in G_{\sig}$. If there is a proper face $\tau \prec \sig$ so that $g \in G_{\tau}$, then $a_{\rho}(g)=1$ for all $\rho \notin \tau(1)$.  In particular, there is $\rho \in \sig(1) \setminus \tau(1)$ with $a_{\rho}(g)=1$. But this fact contradicts $\rho \in g(1)=\sig(1)$, so we get that $g \in G_{\sig}^{\circ}$.
\end{proof}

The case of a simplicial toric variety $X_{\Sig}$ containing a torus factor can be reduced to the previous situation by splitting off a maximal torus factor.
\br Let $X_{\Sig}$ be a simplicial toric variety,
 with $k_{\sig}: V_{\sig}\hookrightarrow X$ the inclusion of the corresponding orbit closure
 for each cone $\sig \in \Sig$. Then the (normalized) mock Hirzebruch class of the orbit closure $V_{\sig}$ (which itself is a simplicial toric variety) can also be computed by the formula:
\be\label{mocko}
\begin{split}\widehat{T}^{(m)}_{y*}(V_{\sig})&= (k_{\sig})^*\left( \prod_{\rho \notin \sig(1)} \frac{[V_{\rho}] \cdot 
\big( 1+y   \cdot e^{-[V_{\rho}](1+y)}\big)}{1-e^{-[V_{\rho}](1+y)}} \right) \cap [V_{\sig}] \\ &=
(k_{\sig})^*\left(\prod_{\rho \notin \sig(1)} \widehat{Q}_y([V_{\rho}] )\right) \cap [V_{\sig}] \: .
\end{split} \ee
Note that (\ref{mock}) corresponds to the case of the trivial cone $\{0\}$. 
\er
For each singular cone $\sig \in \Sig_{\rm sing}$, we set 
\be\label{A} \cA_y({\sig}):= \frac{1}{\mult(\sig)} \cdot \sum_{g \in G_{\sig}^{\circ}}\prod_{\rho \in \sig(1)} 
\frac{ 1+y  \cdot a_{\rho}(g)  \cdot e^{-[V_{\rho}](1+y)}}{1-a_{\rho}(g) \cdot e^{-[V_{\rho}](1+y)}}.\ee
This also agrees by (\ref{Gsig-int}) and (\ref{a-intr}) with the corresponding definition (\ref{A-intro}) given in the introduction. 
 In particular, in the expression of  $\cA_y(\sig)$ only the cohomology classes $[V_{\rho}]$ depend on the fan $\Sig$, but the multiplicity $\mult(\sig)$, the character $a_{\rho}:G_{\sig} \to \bC^*$ and the subset $G_{\sig}^{\circ}$ depend only on rational simplicial cone $\sig$ (and not on the fan $\Sig$ nor the group $G$ from the Cox quotient construction). 
In the above notations, Theorem \ref{Tc} can now be reformulated as follows:
\bt\label{Tcb} The normalized Hirzebruch class of a simplicial toric variety $X=X_{\Sig}$ is computed by the formula:
\be\label{Hirzcb}
\widehat{T}_{y*}(X)=\widehat{T}^{(m)}_{y*}(X) + \sum_{\sig \in \Sig_{\rm sing}}  \cA_y({\sig}) \cdot \left( (k_{\sig})_*\widehat{T}^{(m)}_{y*}(V_{\sig}) \right) \:,
\ee
with $k_{\sig}:V_{\sig} \hookrightarrow X$ denoting the orbit closure inclusion.
\et

\begin{proof} First, by the identification in (\ref{last}), we can rewrite formula (\ref{nHirzc}) as
\be\label{nH}
\widehat{T}_{y*}(X)=\widehat{T}^{(m)}_{y*}(X)+ \left( \sum_{\sig \in \Sig_{\rm sing}} \sum_{g \in G_{\sig}^{\circ}}  \prod_{\rho \in \Sig(1)} \frac{[V_{\rho}] \cdot 
\big( 1+y  \cdot a_{\rho}(g)  \cdot e^{-[V_{\rho}](1+y)}\big)}{1-a_{\rho}(g) \cdot e^{-[V_{\rho}](1+y)}} \right) \cap [X]\:,
\ee
with the mock Hirzebruch class $\widehat{T}^{(m)}_{y*}(X)$ corresponding to the identity element in (\ref{last}). Next, by using Lemma \ref{l10} and the well-known formula from the intersection theory on simplicial toric varieties 
\be\label{int}
\left(\prod_{\rho \in \sig(1)} [V_{\rho}] \right) \cap [X]=\frac{1}{\mult(\sig)} (k_{\sig})_*([V_{\sig}]),
\ee
the contribution of a singular cone $\sig \in \Sig_{\rm sing}$ to $\widehat{T}_{y*}(X)$ can be written as:
\begin{equation*}
\begin{split}
& \left( \sum_{g \in G_{\sig}^{\circ}}  \prod_{\rho \in \Sig(1)} \frac{[V_{\rho}] \cdot 
\big( 1+y  \cdot a_{\rho}(g)  \cdot e^{-[V_{\rho}](1+y)}\big)}{1-a_{\rho}(g) \cdot e^{-[V_{\rho}](1+y)}} \right) \cap [X] \\ 
&= \left(\prod_{\rho \notin \sig(1)} \widehat{Q}_y([V_{\rho}] )\right) \cdot \left[ \sum_{g \in G_{\sig}^{\circ}}  \prod_{\rho \in \sig(1)} 
\frac{ \big( 1+y  \cdot a_{\rho}(g)  \cdot e^{-[V_{\rho}](1+y)}\big)}{1-a_{\rho}(g) \cdot e^{-[V_{\rho}](1+y)}} \right] \cap \left(\big(\prod_{\rho \in \sig(1)} [V_{\rho}] \big)\cap [X] \right) \\
&\overset{(\ref{int})}{=} \cA_y({\sig}) \ \cap \left( \big( \prod_{\rho \notin \sig(1)} \widehat{Q}_y([V_{\rho}] )   \big)   \cap  (k_{\sig})_*([V_{\sig}])\right) \\
&\overset{(*)}{=} \cA_y({\sig}) \ \cap \ (k_{\sig})_*\left((k_{\sig})^*\big(\prod_{\rho \notin \sig(1)} \widehat{Q}_y([V_{\rho}] )\big) \cap [V_{\sig}] \right)\\
&=\cA_y({\sig}) \ \cap \ (k_{\sig})_*\widehat{T}^{(m)}_{y*}(V_{\sig}),
\end{split}
\end{equation*}
where the equality labeled ($*$) follows from the projection formula.
\end{proof}

\br\label{rem_mock} It follows from Theorem \ref{Tcb} that the difference $\widehat{T}_{y*}(X)-\widehat{T}^{(m)}_{y*}(X) $ between the actual and resp. mock Hirzebruch classes 
of the simplicial toric variety $X=X_{\Sig}$ is localized on the singular locus, and the above formula (\ref{Hirzcb}) identifies explicitly the contribution of each singular 
cone to this difference. In particular, since $X$ is smooth in codimension one, and $\widehat{Q}_y(\alpha)=1+\frac{1-y}{2} \al + \cdots$, we get:

$$\widehat{T}_{y*}(X)=[X] + \frac{1-y}{2} \sum_{\rho \in \Sig(1)} [V_{\rho}] + {\rm lower \ order\ homological \ degree \ terms.}$$
In the case of a complete simplicial toric variety, it follows from (\ref{ione}) that the degree zero term of  $\widehat{T}_{y*}(X)$ is $\chi_y(X) \cdot [pt]$. 
It also can be seen from the proof of Theorem \ref{Tc} that for a given singular cone $\sig \in \Sig_{\rm sing}$, the
corresponding mock Hirzebruch class $\widehat{T}^{(m)}_{y*}(V_{\sig})$ of the orbit closure can be regarded as {\it tangential data} for the fixed point sets $W^g$ of the action of $G$ on $W$, whereas the ``coefficient" $\cA_y(\sig)$ in (\ref{Hirzcb}) encodes the {\it normal data} information.
\er

Note that for $y=-1$ one gets by $\cA_{-1}(\sig)=\frac{|G_{\sig}^{\circ}|}{\mult(\sig)}$
 an easy relation between Ehler's formula (\ref{i2}) for the MacPherson Chern class and the mock
 Chern class of a toric variety which is given by
 $$c_*^{(m)}(X_{\Sig}):= \left( \prod_{\rho \in \Sig(1)} (1+ [V_{\rho}])  \right) \cap [X_{\Sig}] 
= \sum_{\sig\in \Sig}\: \frac{1}{\mult(\sig)}\cdot [V(\sig)].$$
 
 Similarly, the top-dimensional contribution of a singular cone $\sig \in \Sig_{\rm sing}$ in  formula
(\ref{Hirzcb}) is given by
\be\label{correction}
\frac{1}{\mult(\sig)} \cdot \sum_{g \in G_{\sig}^{\circ}} \left( \prod_{\rho \in \sig(1)} 
\frac{ 1+y  \cdot a_{\rho}(g)}{1-a_{\rho}(g)} \right) \cdot [V_{\sig}],
\ee
 which for an isolated singularity $V_{\sig}=\{pt\}$ is already the full contribution.
 Specializing further to $y=0$ for a singular cone of smallest codimension
 (so that $G_{\sig}^{\circ}=G_{\sig}\backslash \{id\}$, with $G_{\sig}$ cyclic), we recover the correction factors
 for Todd class formulae of toric varieties appearing in  \cite{BP, BV, GP, P, P2}) also in close relation to (generalized) Dedekind sums.
 \bex\label{exw} Let us now explain the correction term (\ref{correction}) in the case of a $d$-dimensional weighted projective space $X_{\Sig}=\bP(1, \dots, 1,m)$ with weights $(1, \dots, 1,m)$, and $m >1$ and $d>1$, as used in Example \ref{ex-intro} of the Introduction. In this case, the Cox geometric quotient construction is given by the usual action of $\bC^* \simeq G$ on $\bC^{d+1} \setminus \{0\}=W$ (as explained in Example \ref{fake}):
 $$t \cdot (x_0,\dots, x_{d-1}, x_d) = (t x_0,\dots, tx_{d-1}, t^m x_d).$$
 Here, $x_j$ is the coordinate of $W$ corresponding to the ray $\rho_j$ in the associated fan $\Sig$. 
 Then the subset of elements $g$ of $G$ having fixed points is described by $$G_{\Sig}\simeq \mu_m:=\{\lambda \in \bC^*| \ \lambda^m=1\},$$ with fixed point set $W^g=\{x_0=\cdots=x_{d-1}=0\}$ for all $g \in G_{\Sig} \setminus \{id\}$. This fixed point set is mapped by the quotient map to the isolated singular point of $X_{\Sig}$ corresponding to the cone $\sig$ of $\Sig$ spanned by the rays $\rho_j$, $j=0,\dots,d-1$. Since $\lambda$ acts by the usual multiplication on the normal bundle $N_{W^g}$, the contribution of the singular cone $\sigma$ becomes
\be
\frac{1}{\mult(\sig)} \cdot \sum_{g \in G_{\sig}^{\circ}} \left( \prod_{\rho \in \sig(1)} 
\frac{ 1+y  \cdot a_{\rho}(g)}{1-a_{\rho}(g)} \right) \cdot [V_{\sig}]= \frac{1}{m} \sum_{ \lambda^m=1,\\  \lambda \neq 1}  \left( \frac{1+\lambda y}{1-\lambda} \right)^d \cdot [pt]\:.
\ee 
 \eex
 
 \vspace{1.8mm}
 
\br In the above discussion, we used the geometric quotient realization $X_{\Sigma}=W/G$ of a simplicial toric variety (without a torus factor),  as described by the Cox construction. Similar arguments apply also to other such quotient descriptions of $X_{\Sig}$ associated with a stacky fan in the sense of Borison-Chen-Smith \cite{BCS}, with the only modification appearing in the order of stabilizers, e.g., the generic stabilizer need not be trivial. These stabilizers are not intrinsic to the fan $\Sig$, but they depend on the group action, so that multiplicities of cones have to be replaced by orders of the corresponding stabilizer groups. \er

\section{On some characteristic class formulae of Cappell-Shaneson}\label{CS90}

In this final section, we combine the two approaches indicated above for computing homology Hirzebruch classes in order  to derive proofs of several characteristic class formulae, which were first obtained by Cappell and Shaneson, 
see \cite{CS,S}. 

\bd The {\it $T$-class} of a $d$-dimensional simplicial toric variety $X$ is defined as:
\be
T_*(X):=\sum_{k=0}^d 2^{d-k} \cdot td_k(X).
\ee
Note that the $T$-class $T_*(X)$ and the Todd class $td_*(X)$ determine each other.
\ed
The following formula (in terms of $L$-classes) appeared first in Shaneson's $1994$ ICM proceeding paper, see \cite{S}[(5.3)]:
\bt\label{CS1} If $X=X_{\Sig}$ is a $d$-dimensional  simplicial toric variety and for each cone $\sig \in \Sig$ we denote  by $k_{\sig}:V_{\sig} \hookrightarrow X$ the inclusion of the corresponding orbit closure, then 
\be\label{cs1}
T_*(X)=\sum_{\sig \in \Sig} (k_{\sig})_*\widehat{T}_{1*}(V_{\sig}).
\ee
\et
Recall that $L_*(X_{\Sig})=\widehat{T}_{1*}(X_{\Sig})$
for a simplicial projective toric variety or compact toric manifold. Conjecturally, this equality should hold for any complete simplicial toric variety.

The proof of Theorem \ref{CS1} is a direct consequence of the following result:
\bp\label{CS3}
If $X_{\Sig}$ is a simplicial toric variety of dimension $d$, then:
\be\label{cs3}
2^d \cdot td_*(X)=\sum_{\sig \in \Sig} (k_{\sig})_*T_{1*}(V_{\sig}).
\ee
\ep

\begin{proof} By substituting $y=1$ in formula (\ref{7}) of Theorem \ref{t1} we get the following identity:
\be\label{cs5}T_{1*}(X_{\Sig})=\sum_{\sig \in \Sig}\, 2^{\dim(O_{\sig})} \cdot (k_{\sig})_*td_*([\omega_{V_{\sig}}]).\ee
A similar formula holds if we replace $X_{\Sig}$ by an orbit closure $V_{\sig}$. Indeed,  $V_{\sig}$ is itself a toric variety whose fan is build from those cones $\tau \in \Sig$ which have $\sig$ as a face. Recall that if $\sig \preceq \tau$ we have that $O_{\tau} \subset {\bar O}_{\sig}=V_{\sig}$ and, moreover,  $V_{\sig}=\sqcup_{\sig \preceq \tau} O_{\tau}$. Denote by $k_{\tau,\sig}: V_{\tau} \hookrightarrow V_{\sig}$ the induced inclusion. Then formula (\ref{cs5}) applied to $V_{\sig}$ becomes:
\be\label{cs6}
T_{1*}(V_{\sig})=\sum_{\sig \preceq \tau}\, 2^{\dim(O_{\tau})} \cdot (k_{\tau,\sig})_*td_*([\omega_{V_{\tau}}])=\sum_{O_\tau \subset V_{\sig}}\, 2^{\dim(O_{\tau})} \cdot (k_{\tau,\sig})_*td_*([\omega_{V_{\tau}}]).
\ee 
Note that $k_{\sig} \circ k_{\tau,\sig}=k_{\tau}$. So the right-hand side of (\ref{cs3}) can be computed as:
\be\label{cs7}
\begin{split}
\sum_{\sig \in \Sig} (k_{\sig})_*T_{1*}(V_{\sig})&\overset{(\ref{cs6})}{=}\sum_{\sig \in \Sig} \sum_{O_\tau \subset V_{\sig}}\, 2^{\dim(O_{\tau})} \cdot (k_{\tau})_*td_*([\omega_{V_{\tau}}])\\
\end{split}
\ee
On the other hand, for each open orbit $O_{\tau}\subset X$ there are exactly $2^{\rm codim(O_{\tau})}$ open orbits $O_{\sig}$ so that $O_{\tau} \subset {\bar O}_{\sig}=V_{\sig}$. Indeed, this is equivalent to the fact that each simplicial cone $\tau \in \Sig$ has exactly $2^{\dim(\tau)}$ faces $\sig \preceq \tau$ (each of which is generated by a subset of the generators of $\tau$). So the right-hand side of (\ref{cs7}) can be further computed as:
$$ 
\sum_{O_{\tau} \subset X} 2^{\rm codim(O_{\tau})} \cdot 2^{\dim(O_{\tau})} \cdot (k_{\tau})_*td_*([\omega_{V_{\tau}}])=2^d \cdot \sum_{\tau \in \Sig} (k_{\tau})_*td_*([\omega_{V_{\tau}}])\overset{(\ref{T})}{=}2^d \cdot td_*(X),
$$
which completes our proof.
\end{proof}

We can now finish the proof of Theorem \ref{CS1}:
\begin{proof}(of Theorem \ref{CS1})
Let us apply the second homological Adams operation $\Psi_2$ to the formula of Proposition \ref{CS3}. Recall that $\Psi_2$ is defined by multiplying by $2^{-k}$ on $H_{k}(-)\otimes \bQ$. Moreover, since $\Psi_2$ commutes with proper push-down, we get the equality:
\be\label{cs4}
\Psi_2 \big(2^d \cdot td_*(X) \big)=\sum_{\sig \in \Sig} (k_{\sig})_*\Psi_2 \big(T_{1*}(V_{\sig}) \big) = \sum_{\sig \in \Sig} (k_{\sig})_*
\widehat{T}_{1*}(V_{\sig}).
\ee
By definition, the left-hand side of (\ref{cs4}) is just the $T$-class of $X$. Indeed,
$$\Psi_2 \big(2^d \cdot td_*(X) \big)=2^d \cdot  \Psi_2 \big(td_*(X) \big)=2^d \cdot \sum_{k=0}^d 2^{-k} \cdot td_k(X)=\sum_{k=0}^d 2^{d-k} \cdot td_k(X)=T_*(X).$$
\end{proof}

\br 
While Cappell-Shaneson's method of proof of Theorem \ref{CS1} uses mapping formulae for Todd and $L$-classes in the context of resolutions of singularities, our approach above relies on Theorem \ref{t1}, i.e., on the additivity properties of the motivic Hirzebruch classes.
\er

We next describe a result (originally stated in \cite{CS}[Thm.4]) which was used by Cappell and Shaneson for their computation of coefficients of the Ehrhart polynomial of a simplex (see \cite{CS}[Thm.5]). 
Recall that if $X=X_{\Sig}$ is a simplicial toric variety, the mock $L$-classes $L_*^{(m)}$ of $X$ and resp. $V_{\sig}$ are obtained by substituting $y=1$ in the corresponding mock Hirzebruch classes $\widehat{T}^{(m)}_{y*}$ defined in (\ref{mock}) and resp. (\ref{mocko}), that is,
$$
L^{(m)}_{y*}(X):= \left(\prod_{\rho \in \Sig(1)} \frac{[V_{\rho}]}{\tanh([V_{\rho}])} \right) \cap [X] $$
and 
$$
L^{(m)}_{y*}(V_{\sig})= (k_{\sig})^*\left( \prod_{\rho \notin \sig(1)} \frac{[V_{\rho}]}{\tanh([V_{\rho}])} \right)  \cap [V_{\sig}]\:.
$$
\bd(\cite{CS,S}) \ The {\it mock $T$-class} of a  simplicial toric variety $X=X_{\Sig}$ is defined by the formula:
\be\label{mockt}
T_*^{(m)}(X)=\sum_{\sig \in \Sig} (k_{\sig})_*L_*^{(m)}(V_{\sig}).
\ee
Similarly, we define for $\{0\}\neq \sig$ the mock $T$-class of an orbit closure  $V_{\sig}$ of $X$ by:
\be
T_*^{(m)}(V_{\sig})=\sum_{\{\tau| \sig \preceq \tau \}} (k_{\tau,\sig})_*L_*^{(m)}(V_{\tau}),
\ee
where $k_{\tau,\sig}:V_{\tau} \to V_{\sig}$ denotes the map induced by the orbit inclusion $O_\tau \subset V_{\sig}$.
\ed
\br Note that if $X$ is smooth, we get by formula (\ref{cs1}) that $T_*^{(m)}(X)=T_*(X)$.
\er

For each singular cone $\sig \in \Sig_{\rm sing}$, let us now recall the definition of
$ \cA_1({\sig})=\cA_y({\sig})|_{y=1}$ with $\cA_y({\sig})$ defined as in (\ref{A}):
\be\label{al} 
\begin{split}
\cA_1({\sig})&:= \frac{1}{\mult(\sig)} \cdot \sum_{g \in G_{\sig}^{\circ}}\prod_{\rho \in \sig(1)} 
\frac{ 1+   a_{\rho}(g)  \cdot e^{-2[V_{\rho}]}}{1-a_{\rho}(g) \cdot e^{-2[V_{\rho}]}} \\
&=\frac{1}{\mult(\sig)} \cdot \sum_{g \in G_{\sig}^{\circ}}\prod_{\rho \in \sig(1)} 
\frac{ a_{\rho}(g)  \cdot e^{2[V_{\rho}]} +1}{a_{\rho}(g) \cdot e^{2[V_{\rho}]}-1} \:,
\end{split}
\ee
where the second equality follows from the fact that $g \in G_{\sig}^{\circ} \iff g^{-1} \in G_{\sig}^{\circ}$. In fact, the above definition in (\ref{al}) can be extended to all cones of $\Sig$ by setting 
$\cA_y(\{0\})=1$, and $\cA_y({\sig})=0$ for any positive dimensional smooth cone of $\Sig$.

Moreover, it follows by (\ref{Gsig-int}) and (\ref{a-intr}) that in the expression of $\cA_1(\sig)$ (and more generally $\cA_y(\sig)$) only the cohomology classes $[V_{\rho}]$ depend on the fan $\Sig$, while all other quantities depend only on the rational simplicial cone $\sig$ (and not on the fan $\Sig$ nor the group $G$ from the quotient construction) and can be given directly in terms of $\sig$ as already explained in the Introduction. 
In fact, only after this identification our  $\cA_1({\sig})$ above agrees with the corresponding correction factor  $\cA({\sig})$ defined 
just before Theorem 4 in \cite{CS}.

\vspace{2mm}

In the following we need to apply our result from Theorem \ref{Tcb} not just to the given  simplicial toric variety $X$, but also to the orbit closure
$V_{\sig}$, which is again a simplicial toric variety. For $\sig \preceq \tau$,  denote the corresponding correction factors of the closed 
inclusion 
$k_{\tau,\sig}: V_{\tau} \hookrightarrow V_{\sig}$ by $\cA_y({\sig, \tau})$, and similarly for the corresponding group 
$G(\sig,\tau)$, and subset $G(\sig,\tau)^{\circ}$. Then this group and subset  depend only on  the rays
$\{u_{\rho}\}_{\rho\in \tau(1) \setminus \sigma(1)}$. Let $\sig'$ be the cone generated by $\{u_{\rho}\}_{\rho\in \tau(1) \setminus \sigma(1)}$,
so that we get the following commutative diagram of inclusions
\begin{equation*}
\begin{CD}
V_{\tau} @> k_{\tau, \sigma'} >> V_{\sig'}\\
@V k_{\tau, \sig} VV @V k_{\sig'}VV \\
V_{\sig} @> k_{\sig} >> X\:,
   \end{CD}
   \end{equation*}
with $k_{\sig'}\circ k_{\tau, \sigma'} = k_{\tau} = k_{\sig}\circ k_{\tau, \sig}$.
Then $G(\sig,\tau)\simeq G_{\sig'}$, $G(\sig,\tau)^{\circ}\simeq G_{\sig'}^{\circ}$ and
\be\label{A-pullback} \cA_y({\sig, \tau})\simeq (k_{\sig})^*\left( \cA_y({\sig'}) \right)\:.
\ee

We can now prove the following result (stated in \cite{CS}[Thm.4] for a complete simplicial toric variety):
\bt\label{CS2} Let $X=X_{\Sig}$ be a  simplicial toric variety. Then in the above notations we have that:
\be\label{cs2}
\begin{split}
T_*(X)
&= T_*^{(m)}(X)+\sum_{\sig \in \Sig_{\rm sing}} \cA_1(\sig) \cdot (k_{\sig})_*T_*^{(m)}(V_{\sig}) \\
&= \sum_{\sig \in \Sig} \cA_1(\sig) \cdot (k_{\sig})_*T_*^{(m)}(V_{\sig})\:.
\end{split}
\ee
\et

\begin{proof} Recall that by Theorem \ref{CS1} we have that:
$$T_*(X)=\sum_{\sig \in \Sig} (k_{\sig})_*\widehat{T}_{1*}(V_{\sig}).$$
Moreover, by substituting $y=1$ in formula 
(\ref{Hirzcb}) of Theorem \ref{Tcb}, we get:
\be\begin{split} \widehat{T}_{1*}(X)&=L^{(m)}_{*}(X) + \sum_{\sig \in \Sig_{\rm sing}}  \cA_1({\sig}) \cdot  (k_{\sig})_*L^{(m)}_{*}(V_{\sig}) \\ &= 
\sum_{\sig \in \Sig}  \cA_1({\sig}) \cdot  (k_{\sig})_*L^{(m)}_{*}(V_{\sig}) \:,\end{split}\ee
and a similar formula holds for the $L$-class of each orbit closure $V_{\sig}$:
\be\begin{split} \widehat{T}_{1*}(V_{\sig})&=L^{(m)}_{*}(V_{\sig}) + \sum_{\{\tau | \sig \prec \tau\}}  \cA_1({\sig,\tau}) \cdot  (k_{\tau, \sig})_*L^{(m)}_{*}(V_{\tau})\\
&= \sum_{\{\tau | \sig \preceq \tau\}}  \cA_1({\sig, \tau}) \cdot  (k_{\tau, \sig})_*L^{(m)}_{*}(V_{\tau})
\:.\end{split}\ee

The desired formula (\ref{cs2}) follows by combining these facts:
\begin{equation*}
\begin{split}
T_*(X) &= \sum_{\sig \in \Sig} (k_{\sig})_*\widehat{T}_{1*}(V_{\sig}) \\
&= \sum_{\sig \in \Sig} (k_{\sig})_* \left(
\sum_{\{\tau | \sig \preceq \tau\}}  \cA_1({\sig,\tau}) \cdot  (k_{\tau, \sig})_*L^{(m)}_{*}(V_{\tau}) \right) \\
&= \sum_{\sig' \in \Sig} \: \sum_{\sig \in \Sig} (k_{\sig})_* \left(
\quad \sum_{\{\tau | \sig \preceq \tau,\: \tau(1)\setminus \sig(1)\: span\: \sigma'\}} \quad
\cA_1{\sig,\tau}) \cdot  (k_{\tau, \sig})_*L^{(m)}_{*}(V_{\tau}) \right) \\
&\stackrel{(\ref{A-pullback})}{=}  \sum_{\sig' \in \Sig}  \cA_1({\sig'}) \cdot  (k_{\sig'})_*
\left( \sum_{\{\tau | \sig' \preceq \tau\}} (k_{\tau, \sig'})_*L^{(m)}_{*}(V_{\tau})  \right)\\
&= \sum_{\sig' \in \Sig} \cA_1(\sig') \cdot (k_{\sig'})_*T_*^{(m)}(V_{\sig'}) \\
&= T_*^{(m)}(X)+\sum_{\sig' \in \Sig_{\rm sing}} \cA_1(\sig') \cdot (k_{\sig'})_*T_*^{(m)}(V_{\sig'})\:.
\end{split}
\end{equation*}
\end{proof}

\br

While Cappell-Shaneson's approach to Theorem \ref{CS2} uses induction and Atiyah-Singer (resp. Hirzebruch-Zagier) type results for $L$-classes in the {\it local} orbifold description of toric varieties, we use the specialization of Theorem \ref{Tcb} to the value $y=1$, which is based on Edidin-Graham's Lefschetz-Riemann-Roch theorem in the context of the Cox {\it global} geometric quotient construction.

\er

Finally, the following renormalization of Theorem \ref{CS2}  in terms of 
\be\label{alpha} \begin{split} \alpha(\sig):= \Psi_{\frac{1}{2}}(\cA_1(\sig))&=
\frac{1}{\mult(\sig)} \cdot \sum_{g \in G_{\sig}^{\circ}}\prod_{\rho \in \sig(1)} 
\frac{ 1+ a_{\rho}(g)  \cdot e^{-[V_{\rho}]}}{1-a_{\rho}(g) \cdot e^{-[V_{\rho}]}}\\
&= \frac{1}{\mult(\sig)} \cdot \sum_{g \in G_{\sig}^{\circ}}\prod_{\rho \in \sig(1)} 
\coth \left( \pi i \cdot \gamma_{\rho}(g)+\frac{1}{2}[V_{\rho}] \right)
\end{split}
\ee
fits better with the corresponding Euler-MacLaurin formulae (see \cite{CS2}[Thm.2], \cite{S}[Sect.6]), with
$\alpha(\{0\}):=1$ and $\alpha(\sig):=0$ for any other smooth cone $\sig\in \Sig$.
Here the inverse $\Psi_{\frac{1}{2}}$ of the
second (co)homological Adams operation $\Psi_2$  is defined by multiplying by $2^{-k}$ on $H^{k}(-) \otimes \bQ$ and resp. by multiplying by $2^{k}$ on $H_{k}(-) \otimes \bQ$, with 
$$\Psi_{\frac{1}{2}}([V_{\rho}])=\frac{1}{2}[V_{\rho}]\in H^1(X) \otimes \bQ
\quad \text{and} \quad \Psi_{\frac{1}{2}}([X])=2^d[X]\in H_d(X) \otimes \bQ$$
for a simplicial toric variety $X$ of dimension $d$.

\bc Let $X=X_{\Sig}$ be a simplicial toric variety. Then we have:
\be\label{todd-EM}
td_*(X)= \sum_{\sig\in \Sig}\: \alpha(\sig)\cdot \left(
\sum_{\{\tau | \sig \preceq \tau\}}\: \mult(\tau)\prod_{\rho \in \tau(1)} \frac{1}{2}[V_{\rho}]
\prod_{\rho \notin \tau(1)} \frac{\frac{1}{2}[V_{\rho}]}{\tanh(\frac{1}{2}[V_{\rho}])}
\right) \cap [X].
\ee
\ec

\begin{proof}
Let $d$ be the dimension of $X$. Appying $(1/2)^d\cdot \Psi_{\frac{1}{2}}$ to the equality (\ref{cs2}), one gets:
\begin{equation*}\begin{split}
td_*(X)&=\left(\frac{1}{2}\right)^d\cdot \Psi_{\frac{1}{2}}(T(X)) \\
&= \sum_{\sig \in \Sig} \alpha(\sig) \cdot \left(\frac{1}{2}\right)^d\cdot
(k_{\sig})_*\left(  \Psi_{\frac{1}{2}}(T_*^{(m)}(V_{\sig}))\right)\\
&= \sum_{\sig \in \Sig} \alpha(\sig) \cdot \left(\sum_{\{\tau | \sig \preceq \tau\}}\:
\left(\frac{1}{2}\right)^{{\rm codim}(O_{\tau})}\cdot
(k_{\tau})_*\left((k_{\tau})^*(\prod_{\rho \notin \tau(1)} \frac{\frac{1}{2}[V_{\rho}]}{\tanh(\frac{1}{2}[V_{\rho}])})\cap [V_{\tau}]\right)\right)\\
&= \sum_{\sig\in \Sig}\: \alpha(\sig)\cdot \left(
\sum_{\{\tau | \sig \preceq \tau\}}\: \mult(\tau)\prod_{\rho \in \tau(1)} \frac{1}{2}[V_{\rho}]
\prod_{\rho \notin \tau(1)} \frac{\frac{1}{2}[V_{\rho}]}{\tanh(\frac{1}{2}[V_{\rho}])}
\right) \cap [X].
\end{split}
\end{equation*}
Here, the last equality follows from the projection formula for the inclusion
$k_{\tau}: V_{\tau}\hookrightarrow  X=X_{\Sig}$, together with 
the following well-known identity from the intersection theory on simplicial toric varieties: 
$$
\left(\prod_{\rho \in \tau(1)} [V_{\rho}] \right) \cap [X]=\frac{1}{\mult(\tau)} (k_{\tau})_*([V_{\tau}]).
$$
\end{proof}


\end{document}